\documentclass[11pt]{article}
\usepackage{amsthm,amssymb,amsmath,amscd}
\usepackage[all]{xy}
\usepackage{mathrsfs}
\usepackage{multirow}
\usepackage{color}
\usepackage{mathptm}
\usepackage{pgffor}
\usepackage{url}

\objectmargin{1pt}

\newtheorem{Def}{Definition}[section]
\newtheorem{Prop}[Def]{Proposition}
\newtheorem{Theo}[Def]{Theorem}
\newtheorem{Lem}[Def]{Lemma}
\newtheorem{Koro}[Def]{Corollary}
\newtheorem{Bsp}[Def]{Example}

\newtheorem{Rem}[Def]{Remark}

 \DeclareMathOperator{\add}{add}

 \DeclareMathOperator{\Img}{Im}
 
 \DeclareMathOperator{\Ker}{Ker}
 \DeclareMathOperator{\op}{op}
 
 \DeclareMathOperator{\rad}{rad}
 
 \DeclareMathOperator{\soc}{soc}

 \DeclareMathOperator{\tp}{top}

 \DeclareMathOperator{\Hom}{Hom}
 
 \newcommand{\HomP}{\Hom^{\bullet}}

 \newcommand{\otimesP}{\otimes^{\bullet}}
  \DeclareMathOperator{\End}{End}
  
 \DeclareMathOperator{\Ext}{Ext}

\newcommand{\defCategory}[2]{
  \newcommand{#1}{#2\defvariable}}

\newcommand{\defvariable}[2][]{
\if\relax\detokenize{#1}\relax  %if the first arg is empty
\if\relax\detokenize{#2}\relax
    \else  ({#2})  \fi
    \else  ^{{\rm #1}}({#2})  \fi}

\defCategory{\C}{\mathscr{C}}
\defCategory{\K}{\mathscr{K}}
\defCategory{\D}{\mathscr{D}}

 \def\Cb#1{\C[b]{#1}}
 \def\Cf#1{\C[-]{#1}}
 \def\Cz#1{\C[+]{#1}}
 \def\Kb#1{\K[b]{#1}}
 \def\Kf#1{\K[-]{#1}}
 \def\Kz#1{\K[+]{#1}}
 \def\Db#1{\D[b]{#1}}
 \def\Df#1{\D[-]{#1}}
 \def\Dz#1{\D[+]{#1}}

\def\modcat#1{{#1}\mbox{-}{\sf mod}}

\def\pmodcat#1{{#1}\mbox{-}{\sf proj}}

\def\stmodcat#1{{#1}\mbox{-}\underline{\sf mod}}

\newcommand{\lra}{\longrightarrow}
\newcommand{\lla}{\longleftarrow}
\newcommand{\ra}{\rightarrow}
\newcommand{\la}{\leftarrow}
\newcommand{\lraf}[1]{\stackrel{#1}{\lra}}
\newcommand{\llaf}[1]{\stackrel{#1}{\lla}}
\newcommand{\raf}[1]{\stackrel{#1}{\ra}}

\newcommand{\smallvec}[2]{{\displaystyle\scriptsize\left[\begin{matrix}#1\\ #2\end{matrix}\right]}}

\newcommand{\opp}{^{\rm op}}

\newcommand{\cpx}[1]{{#1^{\bullet}}}

\newcommand{\Faf}[1]{#1^{*}}
\newcommand{\nTP}[2]{[{#1}:#2]}

\newcommand{\glueidemef}{ \,{}_e\!\times_{f}}
\newcommand{\glueidem}[2]{ \,{}_{#1}\!\times_{#2}}
\newcommand{\gluesocab}{\, {}_{a}{\diamond}_{b}\,}
\newcommand{\gluesoc}[2]{\, {}_{#1}{\diamond}_{#2}\,}

\newcommand{\stp}[1]{{#1}\mbox{\rm -stp}}

\title{Milnor squares of algebras, I: derived equivalences }
\author{{\sc WEI HU and CHANGCHANG XI$^*$}}
\date{}

\begin{document}
 \maketitle

\renewcommand{\thefootnote}{\alph{footnote}}
\setcounter{footnote}{-1} \footnote{ $^*$ Corresponding author.
Email: xicc@cnu.edu.cn; Fax: 0086 10 68903637.}
\renewcommand{\thefootnote}{\alph{footnote}}
\setcounter{footnote}{-1} 
\renewcommand{\thefootnote}{\alph{footnote}}
\setcounter{footnote}{-1} 

\abstract{Derived equivalences for Artin algebras (and almost $\nu$-stable derived equivalences for finite-dimensional algebras) are constructed from Milnor squares of algebras. Particularly, three operations of gluing vertices, unifying arrows and identifying socle elements on derived equivalent algebras are presented to produce new derived equivalences of the resulting algebras from the given ones. As a byproduct, we construct a series of derived equivalences, showing that derived equivalences may change Frobenius type of algebras in general, though both tilting procedure and almost $\nu$-stable derived equivalences do preserve Frobenius type of algebras.}

\tableofcontents

\section{Introduction}
Pullback algebras (specially, Milnor squares of algebras) appear in many aspects in mathematics. For example, in algebraic $K$-theory, Milnor established a Mayer-Vietoris sequence of $K$-groups for pullback rings (see \cite[Theorem 6.4]{Milnor1971}); in representation theory, Burban and Drozd classified  indecomposable objects of the derived category
of Harish-Chardara modules over $SL({\mathbb R})$
via a special pullback algebra (see \cite{BurbanDrozd}), and Herbara and Prihoda studied infinitely generated projective modules via pullback rings (see \cite{HP2014}); and in homological algebra, Kirkman and Kuzmanovich investigated homological dimensions for pullback
algebras (see \cite{KK}).  One of the important ingredients in these investigations is a characterization of projective modules over a pullback algebra in terms of the ones over its constituent algebras (see \cite[Chapter 2]{Milnor1971}). In the famous Morita theory of derived categories for rings and algebras developed by Rickard (see \cite{Rickard1989, Rickard1989a}), the key notion of tilting complexes involves just a kind of complexes of finitely generated projective modules. This motivates us to consider whether it is possible to get tilting complexes over a pullback algebra through the ones over its constituent algebras. In other words, can we construct derived equivalences by forming pullback algebras?

In this paper, we shall show that under certain conditions, derived equivalences of Artin algebras are preserved by forming pullbacks (see Theorem \ref{TheoMain}). Moreover, if all given derived equivalences are almost $\nu$-stable then so is the induced derived equivalence between pullback algebras (see Corollary \ref{TheoAlmostNuStable}). To apply our result to algebras presented by quivers with relations, we introduce three local operations (gluing vertices, unifying arrows and identifying longest elements) on quivers, so that taking each of them on derived equivalent algebras will produce another derived equivalence of the resulting algebras (see Theorems \ref{TheoremGluing}, \ref{TheoUnifyingArrows} and \ref{theorem-gluesoc}). All of these operations fit well into our framework of constructing derived equivalences for pullback algebras, and can be combined with each other and employed repeatedly. As an application of these techniques, we investigate behaviors of Frobenius parts of derived equivalent algebras and
show that derived equivalences may change Frobenius type of algebras in general.

In Section \ref{sect2}, we fix some notation and recall basic facts needed in later proofs. Particularly, we recall the results on change of rings and on description of projective modules over pullback algebras from \cite{Milnor1971}. Also, we prove some results on images of simple modules under derived equivalences and on tilting complexes and their endomorphism rings.

In Section \ref{sect3},  we first state our main result, Theorem \ref{TheoMain}, which asserts, roughly speaking, that if $A$ is a pullback of homomorphisms $A_1\ra A_0\la A_2$ of Artin algebras with one homomorphism surjective and if $B_i$ is an Artin algebra derived equivalent to $A_i$ for $i=0,1,2$, then there are homomorphisms $B_1\ra B_0\la B_2$ of algebras such that their pullback algebra $B$ is derived equivalent to $A$. After some preparations, we then prove the main result and deduce its corollaries.
Also, we investigate almost $\nu$-stable derived equivalences which induce stable equivalences of Morita type (see \cite{HuXi2010}), and show that, under certain additional conditions, almost $\nu$-stable equivalences between finite-dimensional algebras over an algebraically closed field can be constructed by taking pullback algebras (see Corollary \ref{TheoAlmostNuStable} for precise statement).

In Section \ref{sect4}, we introduce three operations, called gluing vertices, unifying arrows and identifying longest elements, on algebras presented by quivers with relations, and prove that they can produce new derived equivalences from given ones (see Theorems \ref{TheoremGluing}, \ref{TheoUnifyingArrows} and \ref{theorem-gluesoc}). These operations are actually some of effective applications of our main result.

In Section \ref{sect5}, we study, as another application of the main result, the question of whether derived equivalences preserve Frobenius type of algebras. Recall that Frobenius type of algebras means the representation type of their Frobenius parts which have been employed in \cite{HuXi2014-preprint} to lift stable equivalences of Morita type to derived equivalences and in \cite{Martinez-Villa1990a} to reduce the Auslander-Reiten conjecture (or Alperin-Auslander conjecture referred in \cite{Rouquier2006}) on stable equivalences. In this section, we first point out that Frobenius type is  preserved under tilting procedure and almost $\nu$-stable derived equivalences, and then apply our constructions in Section \ref{sect4} to show that derived equivalences may change Frobenius type of algebras in general.

In the second part of this work, we will deal with constructions of stable equivalences of Morita type from pullback algebras.

Acknowledgement. The research work of both authors are partially supported by NSFC. The corresponding author CCX also thanks BNSF
for partial support.
%, while the author W.H. is grateful to the Fundamental Research Funds for the Central Universities for partial support.

\section{Preliminaries\label{sect2}}

In this section, we fix some notation, recall some basic results on derived equivalences and on projective modules over pullback algebras, and then prove a few results concerning derived equivalences and tilting complexes. All results in this section will serve as preparations for the proof of the main result, Theorem \ref{TheoMain}.

\subsection{Derived equivalences}

Let ${\cal C}$ be an additive category.

Given two morphisms
$f:X\rightarrow Y$ and $g:Y\rightarrow Z$ in $\cal C$, the
composite of $f$ with $g$ is written as $fg$, which is a morphism
from $X$ to $Z$. But for two functors $F:\mathcal{C}\rightarrow
\mathcal{D}$ and $G:\mathcal{D}\rightarrow\mathcal{E}$ of
categories, their composition is denoted by $GF$.

For an object $X$ in $\mathcal{C}$, we denote by $\add(X)$ the full
subcategory of $\cal C$ consisting of all direct summands of finite
direct sums of copies of $X$. The object $X$ is said to be \emph{basic} if $X=\bigoplus_{i\in I}X_i$ with $I$ an index set and $X_i$ an indecomposable object for all $i\in I$ such that  $X_i\not\simeq X_j$ for $i\neq j$.

By a complex $\cpx{X}$ over ${\cal C}$ we mean a sequence of morphisms
$d_X^i$ between objects $X^i$ in ${\cal C}: \cdots\ra
X^{i-1}\lraf{d_X^{i-1}}X^i\lraf{d_X^i}X^{i+1}\lraf{d_X^{i+1}}\cdots$,
with $d_X^id_X^{i+1}=0$ for all $i\in\mathbb{Z}$, and write
$\cpx{X}=(X^i, d_X^i)$. The morphism $d_X^i$ is then called the $i$-th differential of $\cpx{X}$.
The complex $\cpx{X}$ is said to be \emph{radical} if each of its differentials is a radical morphism.
By $\cpx{X}[n]$ we denote the $n$-th shift of $\cpx{X}$, that is a complex with the $i$-th term $X^{i+n}$ and differential $(-1)^nd_X^{i+n}$.

We write $\C{\cal C}$ for the category of all complexes over ${\cal C}$, and $\K{\cal C}$ for the homotopy category of $\C{\cal C}$. When ${\cal C}$ is an abelian
category, we write $\D{\cal
C}$ for the derived category of ${\cal C}$. As usual, let $\Cb{\cal C}, \Kb{\cal C}$ and $\Db{\cal C}$ denote the
relevant full subcategories consisting of bounded complexes, respectively; and let
$\Cf{\cal C}, \Kf{\cal C}$ and $\Df{\cal C}$ denote the corresponding full
subcategories consisting of complexes bounded above. Analogously, $\Cz{\cal
C}, \Kz{\cal C}$ and $\Dz{\cal C}$ stand for the corresponding full subcategories
consisting of complexes bounded below, respectively.

Let $\Lambda$ be an Artin algebra over a commutative Artin ring. We denote by $\modcat{\Lambda}$ the
category of finitely generated left $\Lambda$-modules, and by $\pmodcat{\Lambda}$
the full subcategory of $\modcat{\Lambda}$ consisting of finitely
generated projective $\Lambda$-modules. For simplicity, we write
$\C{\Lambda}, \K{\Lambda}$ and $\D{\Lambda}$ for
$\C{\modcat{\Lambda}}, \K{\modcat{\Lambda}}$ and
$\D{\modcat{\Lambda}}$, respectively. Similarly, we have abbreviations
$\Cb{\Lambda}$, $\Kb{\Lambda}$ and $\Db{\Lambda}$. In this paper, $\Db{\Lambda}$ is called the {\em
derived category} of $\Lambda$.

It is well known that the homotopy and derived categories of an Artin algebra (or more generally, a ring) are
triangulated categories. For basic results on triangulated
categories, we refer the reader to the book \cite{Happel1988}.

Two Artin algebras $\Lambda$ and $\Gamma$ are said to be \emph{derived equivalent} if their derived categories are
equivalent as triangulated categories. It follows from Rickard's Morita theory for derived categories \cite{Rickard1989a} that two Artin algebras
$\Lambda$ and $\Gamma$ are derived equivalent if and only if there
is a complex $\cpx{T}$ in $\Kb{\pmodcat{\Lambda}}$ satisfying

\smallskip
(1) $\cpx{T}$ is self-orthogonal, that is, $\Hom_{\Kb{\pmodcat{\Lambda}}}(\cpx{T},\cpx{T}[n])=0$ for all
integers $n\ne 0$;

(2) $\add(\cpx{T})$ generates $\Kb{\pmodcat{\Lambda}}$ as a
triangulated category; and

(3) $\Gamma \simeq \End_{\Kb{\pmodcat{\Lambda}}}(\cpx{T})$.

For more details on derived equivalences, we refer the reader to the papers \cite{Rickard1989, Rickard1989a}, and for some constructions of derived equivalences, we refer the reader to recent papers \cite{HuXi2013, Chen2016}.

A complex $\cpx{T}$ in
$\Kb{\pmodcat{\Lambda}}$ satisfying the above conditions (1) and (2) is
called a {\em tilting complex} over $\Lambda$. It is readily to see that for each tilting complex $\cpx{T}$ there is a basic, radical tilting complex $\cpx{T_0}\in \add(\cpx{T})$ such that $\End_{\Kb{\pmodcat{\Lambda}}}(\cpx{T_0})$ is Morita equivalent to $\End_{\Kb{\pmodcat{\Lambda}}}(\cpx{T}).$

For any derived
equivalence $F:\Db{\Lambda}\ra \Db{\Gamma}$, there is a unique (up
to isomorphism) tilting complex $\cpx{T}$ over $\Lambda$ such that
$F(\cpx{T})$ is isomorphic in $\Db{\Gamma}$ to $\Gamma$. This
complex $\cpx{T}$ is called a tilting complex {\em associated to}
$F$.

Finally, we recall two operations on complexes, which will be used frequently in the paper.

Let $\cpx{X}=(X^i,d_X^i)_{i\in {\mathbb Z}}$ be a complex in
$\C{\Lambda\opp}$ and $\cpx{Y}=(Y^i,d_Y^i)_{i\in {\mathbb Z}}$ a
complex in $\C{\Lambda}$. By $\cpx{X}\otimesP_{\Lambda}\cpx{Y}$ we
mean the total complex of the double complex with $(i,j)$-term
$X^i\otimes_{\Lambda}Y^j$. That is, the $n$-th term of the complex
$\cpx{X}\otimesP_{\Lambda}\cpx{Y}$ is $\bigoplus_{p+q=n}
X^p\otimes_{\Lambda}Y^q =\bigoplus_{q\in {\mathbb
Z}}X^{n-q}\otimes_{\Lambda}Y^{q}$, and the $n$-th differential is given by
$x\otimes y\mapsto x\otimes
(y)d_Y^{q}+(-1)^q(x)d_X^{n-q}\otimes y$ for $x\in X^{n-q}$
and $y\in Y^q$.

Let $\cpx{X}$ and $\cpx{Y}$ be two complexes in $\C{\Lambda}$. By
$\HomP_{\Lambda}(\cpx{X},\cpx{Y})$ we denote the total complex of
the double complex with $(i,j)$-term $\Hom_{\Lambda}(X^{-i},
Y^{j})$. Thus the $n$-th term of the complex
$\HomP_{\Lambda}(\cpx{X},\cpx{Y})$ is $\prod_{p\in {\mathbb
Z}}\Hom_{\Lambda}(X^p,Y^{n+p})$, and the $n$-th differential is given by
$(\alpha^p)_{p\in {\mathbb Z}}\mapsto (\alpha^pd_Y^{n+p}-(-1)^nd_X^p\alpha^{p+1})_{p\in {\mathbb Z}}$
for $\alpha^p\in \Hom_{\Lambda}(X^p,Y^{n+p}).$

\subsection{Complexes under change of rings}

Let $f: \Lambda\ra
\Gamma$ be a homomorphism of $R$-algebras, where $R$ is a commutative ring with identity.
Then every $\Gamma$-module $U$ can be viewed as a $\Lambda$-module by defining
$a\cdot u:=(a)f u$ for all $a\in\Lambda$ and $u\in U$. Thus we get the so-called restriction functor ${}_{\Lambda}(-): \modcat{\Gamma}\ra\modcat{\Lambda}$.  Moreover, there is an adjoint pair  $(\Gamma\otimes_{\Lambda}-, {}_{\Lambda}(-))$ of functors whose unit is the canonical homomorphism of $\Lambda$-modules:
$$f^*: X\lra {}_{\Lambda}\Gamma\otimes_{\Lambda}X, \quad x\mapsto 1\otimes x \; \mbox{ for } \; x\in X.$$

The following lemma tells us about change of projective modules.

\begin{Lem}
Let $f: \Lambda\ra \Gamma$ be a surjective homomorphism of Artin algebras. Then
$\Gamma\otimes_{\Lambda}-$ gives a one-one correspondence between
the set of isomorphism classes of indecomposable projective
$\Lambda$-modules $X$ with $\Gamma\otimes_{\Lambda}X\neq 0$ and the
set of isomorphism classes of indecomposable projective
$\Gamma$-modules.
\label{LemSurjAlgHom}
\end{Lem}

The following lemma is standard for change of rings .

\begin{Lem}
 Let $f: \Lambda\ra \Gamma$ be a homomorphism of Artin algebras,
and let $X$ be a $\Lambda$-module and $U$ a $\Gamma$-module.

 $(1)$ If $f$ is surjective, then so is $f^*: X\ra \Gamma\otimes_{\Lambda}X$.

 $(2)$ There is a natural isomorphism $\Hom_{\Gamma}(\Gamma\otimes_{\Lambda}X, U)\ra \Hom_{\Lambda}(X, U)$ sending $g$ to $f^*g$.
 \label{LemAlgebraHommorphism}
\end{Lem}

\medskip
Using Lemma \ref{LemAlgebraHommorphism}, we can extend results
on modules to complexes. Let $f: \Lambda\ra \Gamma$ be a homomorphism of Artin algebras. Then we have a functor $\Gamma\otimesP_{\Lambda}-: \C{\Lambda}\ra \C{\Gamma}$, which has the restriction functor as its right adjoint functor. So, the unit of this adjoint pair of functors provides a natural chain map $f^*: \cpx{X}\ra
\Gamma\otimesP_{\Lambda}\cpx{X}$ for $\cpx{X}\in \C{\Lambda}$. More precisely, $f^*$ is defined by $f^i: X^i\ra
\Gamma\otimes_{\Lambda}X^i$ for all integers $i$.  As in the case of modules, we have the following lemma for
complexes. Its proof is just a consequence of the universal properties of units of adjoint functors.

\begin{Lem}
  Let $f:\Lambda\ra \Gamma$ be a homomorphism between
  Artin algebras $\Lambda$ and $\Gamma$. Then, for any $\cpx{X}\in \C{\Lambda}$ and $\cpx{U}\in\C{\Gamma}$, we have the following:

  $(1)$ The morphism $\Hom_{\C{\Gamma}}(\Gamma\otimesP_{\Lambda}\cpx{X}, \cpx{U})\ra\Hom_{\C{\Lambda}}(\cpx{X}, \cpx{U})$
     sending $\cpx{h}$ to $f^*\cpx{h}$ is a natural isomorphism.

  $(2)$ The morphism $\Hom_{\K{\Gamma}}(\Gamma\otimesP_{\Lambda}\cpx{X}, \cpx{U})\ra\Hom_{\K{\Lambda}}(\cpx{X},
  \cpx{U})$ sending $\cpx{h}$ to $f^*\cpx{h}$ is a natural
  isomorphism.

  $(3)$ If $\,\cpx{U}\simeq
  \Gamma\otimesP_{\Lambda}\cpx{X}$ in $\C{\Gamma}$, then, for each
  epimorphism $\cpx{g}: \cpx{X}\ra {}_{\Lambda}\cpx{U}$ in $\C{\Lambda}$, there is
  an isomorphism $\cpx{h}: \Gamma\otimesP_{\Lambda}\cpx{X}\ra
  \cpx{U}$ in $\C{\Gamma}$ such that $\cpx{g}=f^*\cpx{h}$.
  \label{LemmaAlghomForComplex}
\end{Lem}

Let $\cpx{X}\in \C{\Lambda}$, $\cpx{U}\in \C{\Gamma}$ and $\cpx{g}: \cpx{X}\ra {}_{\Lambda}\cpx{U}$ be a chain map in $\C{\Lambda}$. If, for each morphism $\cpx{\alpha}: \cpx{X}\ra \cpx{X}$ in $\K{\Lambda}$, there is a unique morphism $\cpx{\beta}: \cpx{U}\ra\cpx{U}$ in $\K{\Gamma}$ such that $\cpx{g}\cpx{\beta}=\cpx{\alpha}\cpx{g}$ in $\K{\Lambda}$, then the map
$$\End_{\K{\Lambda}}(\cpx{X})\lra \End_{\K{\Gamma}}(\cpx{U})$$
sending $\cpx{\alpha}$ to $\cpx{\beta}$ is a homomorphism of algebras, which is called the \emph{algebra homomorphism determined by} $\cpx{g}$.
According to Lemma \ref{LemmaAlghomForComplex}(2), the morphism
$f^*: \cpx{X}\ra \Gamma\otimesP_{\Lambda}\cpx{X}$
determines an algebra homomorphism:
$$ \End_{\K{\Lambda}}(\cpx{X})\lra\End_{\K{\Gamma}}(\Gamma\otimesP_{\Lambda}\cpx{X}).$$
By universal property of units of adjoint functors, we know that the above homomorphism of algebras is actually given by $\alpha\mapsto \Gamma\otimesP_{\Lambda}\alpha$ for $\alpha\in \End_{\K{\Lambda}}(\cpx{X})$.

\subsection{Simple modules under derived equivalences}

Let $\Lambda$ be an  Artin algebra and $Y$ be an
indecomposable $\Lambda$-module. For each $\Lambda$-module $X$, we
decompose $X$ into a direct sum of indecomposable modules, say $X =
\bigoplus_{i=1}^nX_i$, and let $[X: Y]$ be the multiplicity of $Y$ as a
direct summand of $X$, that is, the number of those $X_j$ with $X_j\simeq Y$. Note that $[X:Y]$ is independent of the
choice of the decomposition of $X$. For a bounded complex $\cpx{X}$
over $\modcat{\Lambda}$, we define
$$\nTP{\cpx{X}}{Y}:=\displaystyle{\sum_{i\in\mathbb{Z}}[X^i: Y]}.$$
Note that $\nTP{\cpx{X}}{Y}$ is well defined in
$\Cb{\Lambda}$ by the Krull-Remak-Schmidt Theorem. In the following, we denote by $S_P$ the top of a projective module $P$.

\begin{Lem}
  Let $\cpx{T}$ be a basic,
  radical tilting complex over $\Lambda$, and let
  $\Gamma:=\End_{\Kb{\Lambda}}(\cpx{T})$. Suppose that $F:\Db{\Lambda}\ra\Db{\Gamma}$ is a
  derived equivalence induced by $\cpx{T}$ and that $P$ is an
  indecomposable projective $\Lambda$-module. Then $F(S_P)$ is
  isomorphic in $\Db{\Gamma}$ to $S[n]$ for some simple
  $\Gamma$-module $S$ and some integer $n$ if and only if
  $\nTP{\cpx{T}}{P}=1$.
 \label{LemsimpleTosimple}
\end{Lem}

{\it Proof.}
Suppose $\nTP{\cpx{T}}{P}=1$. Then there is some integer
$n$ such that $\Hom_{\K{\Lambda}}(\cpx{T}, S_P[i])= 0$ for all
$i\neq -n$. Hence $F(S_P[-n])$ is isomorphic in
$\Db{\Gamma}$ to a $\Gamma$-module $X$. Now we prove that $X$ is
simple. Since $\nTP{\cpx{T}}{P}=1$, there is only one indecomposable
direct summand $\cpx{T_P}$ of $\cpx{T}$ such that $P$ occurs in
$\cpx{T_P}$. Let $\bar{P}$ be the indecomposable projective
$\Gamma$-module $F(\cpx{T_P})$. Then $\Hom_{\Db{\Lambda}}(\cpx{T},
S_P[-n])\simeq \Hom_{\Db{\Lambda}}(\cpx{T_P}, S_P[-n])$,
or equivalently $\Hom_{\Gamma}(\Gamma, X)\simeq
\Hom_{\Gamma}(\bar{P}, X)$. This means that $X$ only contains
composition factors isomorphic to $S_{\bar{P}}$. Moreover,
$\End_{\Gamma}(X)\simeq\End_{\Lambda}(S_P)$ is a division
algebra. Hence $X$ must be simple. Note that in the foregoing proof we only need $\cpx{T}$ to be a tilting complex.

Conversely, suppose that
$F(S_P[k])$ is isomorphic to a simple $\Gamma$-module $S$ for
some integer $k$. Then, by
assumption, $\Gamma$ is a basic algebra and $S$ is a $1$-dimensional module over
$D:=\End_{\Gamma}(S)$. Thus  $\Hom_{\K{\Lambda}}(\cpx{T},
S_P[i])$ is zero for all $i\neq k$, and $1$-dimensional over
$D$ for $i=k$. Since $\cpx{T}$ is a radical complex,
$$\Hom_{\K{\Lambda}}(\cpx{T}, S_P[i])\simeq \Hom_{\Lambda}(T^{-i},
S_P)$$ for all integers $i$. This implies that the indecomposable
projective module $P$ occurs in $\cpx{T}$ only in degree $-k$ with  the
multiplicity $1$. Hence $\nTP{\cpx{T}}{P}=1$.
$\square$

\medskip
As an immediate consequence of the proof of Lemma \ref{LemsimpleTosimple}, we get the following corollary for tilting modules.

\begin{Koro} If $T=P\oplus P'$ is a basic titling $\Lambda$-module, where $P$ is projective and $P'$ has a minimal projective resolution $\cpx{Q}=(Q^i,d^i)_{i\le 0}$ such that each indecomposable direct summand of $P$ does not appear in $\bigoplus_{i\ge 0}Q^{-i}$, then there exists a derived equivalence $F: \Db{\Lambda}\ra\Db{\End_{\Lambda}(T)}$ such that $F(S)$ is isomorphic to a simple $\End_{\Lambda}(T)$-module for all simple modules $S\in\add(S_P)$.
\end{Koro}

The following lemma is very useful in our later proofs.

\begin{Lem}
Let $\{U_1, \cdots, U_s, V_1,\cdots,
V_r\}$ be a complete set of pairwise non-isomorphic
   indecomposable projective $\Lambda$-modules and let $U:=\bigoplus_{i=1}^sU_i$. Suppose that
   $\cpx{T}$ is a basic, radical tilting complex over $\Lambda$ with $\nTP{\cpx{T}}{V_i}=1$ for all $1\le i\le r$. Then $\cpx{T}$ can be written as a direct sum (in $\Kb{\Lambda}$)
   $$\cpx{T}\simeq \cpx{U}\oplus
\cpx{V_1}\oplus\cdots\oplus\cpx{V_r}$$ of complexes $\cpx{U}$ and $\cpx{V_i}$ with $1\le i\le r$, satisfying the following properties:

   $({\rm a})$ $\nTP{\cpx{V_i}}{V_j}=1$ for $i=j$, and zero otherwise. Moreover, all $\cpx{V_i}$ are indecomposable complexes.

   $({\rm b})$ $\cpx{U}\in\Kb{\add(U)}$, and $\add(\cpx{U})$ generates $\Kb{\add(U)}$ as a triangulated category.
\label{LemFormofTiltingComp}
\end{Lem}

{\it Proof.}
Let $\Gamma:=\End_{\K{\Lambda}}(\cpx{T})$ and
$F:\Db{\Lambda}\ra\Db{\Gamma}$ be a derived equivalence induced by
the tilting complex $\cpx{T}$. By Lemma \ref{LemsimpleTosimple},
there are pairwise non-isomorphic indecomposable projective
$\Gamma$-modules $\bar{V}_1, \cdots, \bar{V}_r$ such that
$F(\tp(V_i))\simeq \tp(\bar{V_i})[n_i]$ for some $n_i$ with $1\le i\le r$. Let $\bar{U}_1, \cdots, \bar{U}_s$ be
indecomposable projective $\Gamma$-modules such that $\{\bar{U}_1,
\cdots, \bar{U}_s, \bar{V}_1, \cdots, \bar{V}_r\}$ is a complete set
of pairwise non-isomorphic indecomposable projective
$\Gamma$-modules and set $\bar{U}:=\bigoplus_{i=1}^s\bar{U}_i$.
Since $\cpx{T}$ is a basic tilting complex, $\Gamma$ is a
basic algebra, and therefore
$$_{\Gamma}\Gamma\simeq \bar{U}\oplus
\bar{V}_1\oplus\cdots\oplus\bar{V}_r.$$ By definition,
$F(\cpx{T})\simeq {}_{\Gamma}\Gamma$. Now, let $\cpx{U}$ be a direct
summand of $\cpx{T}$ such that $F(\cpx{U})\simeq \bar{U}$ and let
$\cpx{V_i}$ be a direct summand of  $\cpx{T}$ such that
$F(\cpx{V_i})\simeq \bar{V_i}$ for $1\le i\le r.$ Then $F(\cpx{U}\oplus
\cpx{V_1}\oplus\cdots\oplus\cpx{V_r}) \simeq {}_{\Gamma}\Gamma\simeq F(\cpx{T})$, and consequently
 $$\cpx{T}\simeq \cpx{U}\oplus
\cpx{V_1}\oplus\cdots\oplus\cpx{V_r} \; \mbox{ in } \; \Db{\Lambda}.$$ This implies $\cpx{T}\simeq \cpx{U}\oplus
\cpx{V_1}\oplus\cdots\oplus\cpx{V_r} \; \mbox{ in } \; \Kb{\Lambda}.$ Now we have
$$\Hom_{\Kb{\Lambda}}(\cpx{V_i}, \tp(V_j)[k])\simeq \Hom_{\Db{\Lambda}}(\cpx{V_i}, \tp(V_j)[k])\simeq \Hom_{\Db{\Gamma}}(\bar{V_i}, \tp(\bar{V_j})[k+n_j])=0$$
whenever $i\neq j$ or $k\neq -n_j$. By assumption,
$\nTP{\cpx{T}}{V_i}=1$ for  $1\le i\le r$. This forces that the
projective module $V_i$ only occurs in the $(-n_i)$-th degree of
$\cpx{V_i}$.

Now, it is easy to see that all complexes $\cpx{V_i}$ can be chosen to be indecomposable. This proves (a).

By (a) and our assumption $\nTP{\cpx{T}}{V_i}=1$ for all $i$,
the complex $\cpx{U}$ is clearly in $\Kb{\add(U)}$. Now we show that
$F$ induces a triangle equivalence between $\Kb{\add(U)}$ and
$\Kb{\add(\bar{U})}$. In fact, a complex $\cpx{P}$ from
$\Kb{\pmodcat{\Lambda}}$ lies in $\Kb{\add(U)}$ if and only if
$\Hom_{\Db{\Lambda}}(\cpx{P}, \tp(V_i)[k])=0$ for all
$1\le i\le r$  and all $k\in\mathbb{Z}$. However, this is
equivalent to $\Hom_{\Db{\Gamma}}(F(\cpx{P}),
\tp(\bar{V_i})[k+n_i])=0$ for all $1\le i\le r$ and all
$k\in\mathbb{Z}$, that is, $F(\cpx{P})$ belongs to
$\Kb{\add(\bar{U})}$. Hence $F$ induces a triangle equivalence
between $\Kb{\add(U)}$ and $\Kb{\add(\bar{U})}$. Since
$\add(\bar{U})$ generates $\Kb{\add(\bar{U})}$ as a triangulated
category, $\add(\cpx{U})$ generates $\Kb{\add(U)}$ as a
triangulated category. This proves (b).
$\square$

\subsection{Projective modules over Milnor squares of algebras}\label{SectionPullback}

Let $A_0$, $A_1$ and $A_2$ be rings with identity. Given two homomorphisms $\pi_i: A_i\ra A_0$ of rings,
%$$\xymatrix{
%& *+[r]{A_1}\ar[d]^{\pi_1}\\
%A_2\ar[r]^{\pi_2} & *+[r]{A_0}, }$$
the \emph{pullback ring} $A$ of $\pi_1$ and $\pi_2$ is defined by $A:=\{(x,y)\in A_1\times
A_2\mid (x)\pi_1=(y)\pi_2\}$. Transparently, there is a
commutative diagram of ring homomorphisms
$$\xymatrix@M=0.6mm{
  A\ar[r]^{\lambda_1}\ar[d]_{\lambda_2} & *+[r]{A_1}\ar[d]^{\pi_1}\\
  *+[r]{A_2} \ar[r]^{\pi_2} & *+[r]{A_0}
}$$
where $\lambda_i$ is the canonical projections from $A$ to $A_i$ for $i=1,2.$ The above pullback diagram has a universal property:
For any ring homomorphisms $i_1: B\ra A_1$ and $i_2:B\ra A_2$ with $i_1\pi_1=i_2\pi_2$, there is a unique ring homomorphism
$\theta:B\ra A$ such that $\theta \lambda_j=i_j$ for $j=1,2$. Note that if $\pi_1$ is
surjective then so is $\lambda_2$.

If, additionally, one of $\pi_1$ and $\pi_2$ is  surjective, then the above square is called a \emph{ Milnor square} of rings (see \cite{Milnor1971}).
For a Milnor square of rings, there is a nice description of projective
$A$-modules via projective $A_i$-modules in \cite{Milnor1971}. Let us recall it right now.

Given a projective $A_1$-module $X_1$, a projective $A_2$-module
$X_2$ and an isomorphism
$h: A_0\otimes_{A_1}X_1\ra A_0\otimes_{A_2}X_2$ of
$A_0$-modules, the \emph{Milnor patching} of the triple $(X_1,X_2,h)$ is defined by
$$M(X_1, X_2, h):=\big\{(x_1, x_2)\in X_1\oplus
X_2\mid (x_1)\Faf{\pi_1}h=(x_2)\Faf{\pi_2}\big\}=\{(x_1,x_2)\in X_1\oplus X_2\mid (1\otimes x_1)h=1\otimes x_2\}.$$ Let $p_i: M(X_1, X_2, h)\ra X_i$ be the canonical projection. Note that
$M(X_1, X_2, h)$ has an $A$-module structure: For $a\in A,$
$$a\cdot (x_1, x_2)=\big((a)\lambda_1\cdot x_1,\, (a)\lambda_2\cdot x_2\big)\; \mbox{ for }\; x_1\in X_1, x_2\in X_2.$$

Now, we state the following description of projective $A$-modules given in \cite[Chapter 2]{Milnor1971}.
\begin{Lem}
Suppose that $\pi_1$ is surjective, $X_i$ is a projective $A_i$-module
for $i=1,2$, and
$h: A_0\otimes_{A_1}X_1\ra A_0\otimes_{A_2}X_2$ is an isomorphism of
$A_0$-modules. Then we have the following:

$(1)$ The module $M(X_1, X_2, h)$ is a projective $A$-module.
Furthermore, if, in addition, $X_1$ and $X_2$ are finitely generated over $A_1$
and $A_2$, respectively, then $M(X_1, X_2, h)$ is finitely generated
over $A$.

$(2)$ Every projective $A$-module is isomorphic to $M(X_1, X_2, h)$
for some suitably chosen $X_1, X_2$ and $h$.

$(3)$ For $i\in\{1,2\}$, there is a natural isomorphism
$$\mu_i: A_i\otimes_AM(X_1, X_2, h)\lra X_i$$
sending $a_i\otimes (x_1, x_2)$ to $a_ix_i$, and the canonical
projection $p_i: M(X_1, X_2, h)\ra X_i$ is equal to
$\lambda_i^*\mu_i$.

$(4)$ There is an exact sequence of $A$-modules:
$$\xymatrix@M=1.5mm{0\ar[r]&  M(X_1, X_2, h)\ar[r]^(.55){[p_1, p_2]} & X_1\oplus
X_2\ar[r]^{\smallvec{\Faf{\pi_1}h}{\Faf{-\pi_2}}} &A_0\otimes_{A_2}X_2\ar[r] & 0}.$$

\label{LemmaPullbackProj}
\end{Lem}

{\it Proof.}
The statements (1), (2) and (3) are just \cite[Theorems
2.1, 2.2, and 2.3, p. 20]{Milnor1971}. The statement (4)
follows easily from the definition of $M(X_1, X_2, h)$ and the fact
that $\pi_1^*$ is surjective.
$\square$

\medskip
For the rest of this section, we shall assume that $A_0, A_1$ and
$A_2$ are Artin algebras and that
$\pi_1$ is surjective. Thus we have an exact sequence of $A$-bimodules:
$$(*)\quad \xymatrix@M=1mm{
0\ar[r] & A \ar[r]^(.35){[\lambda_1, \lambda_2]} & A_1\oplus A_2 \ar[r]^{\smallvec{\pi_1}{-\pi_2}} & A_0\ar[r] & 0.
}$$

In the following, we shall give a partition of indecomposable projective
$A$-modules.

Let $P_1$ be a direct sum
of all non-isomorphic indecomposable projective $A_1$-modules $X$
such that $A_0\otimes_{A_1}X= 0$, and let $Q_1$ be a direct sum of
all non-isomorphic indecomposable projective $A_1$-modules $Y$ such that
$A_0\otimes_{A_1}Y\neq 0$. Thus $\pmodcat{A_1}=\add(P_1\oplus Q_1)$. Similarly, we define
projective $A_2$-modules  $P_2$ and $Q_2$, and get $\pmodcat{A_2}=\add(P_2\oplus Q_2)$.

Since $\pi_1$ is surjective, $\lambda_2$ is also
surjective. Therefore, if $X$ is an indecomposable projective
$A$-module with $A_2\otimes_AX\neq 0$, then $A_2\otimes_AX$ is an indecomposable projective $A_2$-module by Lemma \ref{LemSurjAlgHom}. Hence, for an
indecomposable projective $A$-module $X$, only the following
three cases occur:

\medskip
$\bullet$ Case 1:  $A_2\otimes_AX=0.$

$\bullet$ Case 2: $0\ne A_2\otimes_AX\in\add(P_2)$.

$\bullet$ Case 3: $0\ne A_2\otimes_AX\in\add(Q_2)$.

\medskip
According to the three cases, we have a partition of indecomposable
projective $A$-modules: For $1\le i\le 3$, let $F_i$ be the direct sum of all
non-isomorphic indecomposable projective $A$-modules $X$ corresponding to Case $i$.
Then $\pmodcat{A} = \add(F_1\oplus F_2\oplus F_3)$.

\begin{Lem}
With the above notation, we have the following:

$(1)$ The functor $A_i\otimes_A-$ and the restriction functor
$_A(-)$ induce mutually inverse equivalences between $\add(F_i)$ and
$\add(P_i)$ for $i=1,2.$

$(2)$ Let $i\in \{1,2\}$ and $X\in \add(P_i)$. Then the natural map
$\Hom_A({}_AX, A)\ra \Hom_{A_i}(X, A_i),$ sending $\alpha$ to
$\alpha\lambda_i$, is an isomorphism of right $A$-modules.

$(3)$ Let $i\in \{1,2\}$ and $X\in \add(P_i)$. If
$\add({}_{A_i}X)=\add(\nu_{A_i}X)$, then $\add({}_AX)=\add(\nu_AX)$, where $\nu_A$ is the Nakayama functor $D\Hom_A(-,{}_AA)$ of $A$.
 \label{LemAddP}
\end{Lem}

{\it Proof.}
(1) We prove the case $i=1$. For each $X$ in $\add(F_1)$, we have
$A_2\otimes_AX=0$, and therefore
$A_0\otimes_{A_1}A_1\otimes_AX\simeq
A_0\otimes_{A_2}A_2\otimes_AX=0$ and $A_1\otimes_AX\in\add(P_1)$.
Thus $X\simeq M(A_1\otimes_AX, 0, 0)$ and the map $\lambda_1^*:
X\ra A_1\otimes_AX$ is a bijection by the definition of
$M(A_1\otimes_AX, 0, 0)$. It follows from  the statements (1) and (3) in Lemma
\ref{LemmaPullbackProj} that, for $X$ and $Y$ in
$\add(F_1)$, the functor $A_1\otimes_A-$ induces an isomorphism
from $\Hom_A(X, Y)$ to $\Hom_{A_1}(A_1\otimes_AX,
A_1\otimes_AY)$. Moreover, for each $U$ in $\add(P_1)$, the module
$M(U, 0, 0)$ is a projective module in $\add(F_1)$ such that
$A_1\otimes_AM(U, 0, 0)\simeq U$. This shows that the functor $A_1\otimes_A-:
\add(F_1)\ra \add(P_1)$ is an equivalence. Clearly, the restriction
functor $_A(-)$ is right adjoint to $A_1\otimes_A-$ by Lemma
\ref{LemAlgebraHommorphism}(2), and therefore a quasi-inverse of
$A_1\otimes_A-$. This proves (1) for $i=1$. The case $i=2$ can be shown similarly.

(2) Assume both $i=1$ and $X\in\add(P_1)$. By Lemma
\ref{LemAlgebraHommorphism}, $\Hom_A(_{A}X, A_i)\simeq
\Hom_{A_i}(A_i\otimes_AX, A_i)$ for $0\le i\le 2$. Since
$X\in\add(P_1)$,  we have $_AX\in\add(F_1)$ by (1), and consequently
$A_2\otimes_AX=0$ and $A_0\otimes_AX\simeq
A_0\otimes_{A_2}A_2\otimes_AX=0$. Therefore  $\Hom_A({}_AX, A_0)=0=\Hom_A({}_AX, A_2)$. Applying $\Hom_A({}_AX, -)$ to
the exact sequence ($*$)
of $A$-bimodules, we get an isomorphism of right $A$-modules:
$$\Hom_A({}_AX, A)\lra \Hom_A({}_AX, {}_AA_1),$$
which sends $\alpha$ to $\alpha\lambda_1$. For the case $i=2$, a proof can be demonstrated similarly.

(3) Without loss of generality, we can assume that the module $X$ is basic. Then it follows from $\add({}_{A_i}X)=\add(\nu_{A_i}X)$ that $\nu_{A_i}X\simeq X$. This together with (2) implies the following isomorphisms:
$$\nu_AX=D\Hom_{A}({}_AX, A)\simeq D(\Hom_{A_i}(X, A_i)_{A})={}_A(\nu_{A_i}X)\simeq {}_AX.$$
Thus (3) follows.
$\square$

The next lemma describes indecomposable projective $A$-modules in $\add(F_3)$.

\begin{Lem}
$(1)$ For each indecomposable $A_2$-module $V$ in $\add(Q_2)$, there
is an $A_1$-module $W$ $($unique up to isomorphism$)$ in $\add(Q_1)$
with an isomorphism $h: A_0\otimes_{A_1}W\ra A_0\otimes_{A_2}V$ such
that  $M(W, V, h)$ is an indecomposable projective $A$-module in
$\add(F_3)$.

$(2)$ Let $\{V_1, \cdots, V_s\}$ be a complete set of pairwise
non-isomorphic indecomposable projective $A_2$-modules in
$\add(Q_2)$, and let $W_i\in \add(Q_1)$ be the projective $A_1$-module determined by $V_i$ in \emph{(1)} for $1\le i\le s$. Then $\{M(W_i, V_i,
h_i) \mid 1\le i\le s\}$ is a complete set of pairwise non-isomorphic
indecomposable projective $A$-modules in $\add(F_3)$.
\label{LemAddQ}
\end{Lem}

{\it Proof. }  (1) Since $\pi_1$ is surjective, it follows from Lemma
\ref{LemSurjAlgHom} that there is an $A_1$-module $W$ (unique up to
isomorphism) and an isomorphism $h:
A_0\otimes_{A_1}W\ra A_0\otimes_{A_2}V$. We need to show that $M(W,
V, h)$ is in $\add(F_3)$. Let $X$ be an indecomposable direct
summand of $M(W, V, h)$. Then there are two possibilities: $A_2\otimes_AX\ne 0$ or $A_2\otimes_AX = 0$. If $A_2\otimes_AX\neq 0$, then $A_2\otimes_AX$
is a direct summand of $V$. Since $V$ is indecomposable, we have
$A_2\otimes_AX\simeq V$. By definition, $X\in\add(F_3)$. Now, we exclude the case
$A_2\otimes_AX=0$. If this happens, then $A_1\otimes_AX\ne 0$. Otherwise $X\simeq M(A_1\otimes_AX,
A_2\otimes_AX, g)=0$. So $A_1\otimes_AX$ is a nonzero direct
summand of $W$. However, by definition, $X\in\add(F_1)$. It follows from
Lemma \ref{LemAddP}(1) that the module $A_1\otimes_AX$ lies in $\add(P_1)$.
This is a contradiction. Thus $M(W, V, h)\in \add(F_3)$. Since
$A_2\otimes_AM(W, V, h)\simeq V$ is indecomposable, the module $M(W, V, h)$ is indecomposable by Lemma
\ref{LemSurjAlgHom}.

(2) It follows from (1) that $M(W_i, V_i, h_i)\in\add(F_3)$ is
indecomposable for all $1\le i\le s$. Now, let $X$ be an indecomposable
$A$-module in $\add(F_3)$. Then the $A_2$-module $A_2\otimes_AX$ is
indecomposable since $\lambda_2$ is surjective. Thus, there is some
$V_i$ such that $A_2\otimes_AX\simeq V_i\simeq A_2\otimes_AM(W_i,
V_i, h_i)$. By Lemma \ref{LemSurjAlgHom}, we have $X\simeq M(W_i,
V_i, h_i)$. This finishes the proof. $\square$

\medskip
Finally, we extend previous facts on modules to the ones on complexes of modules.

Given a complex $\cpx{X_1}$ in $\Cb{\pmodcat{A_1}}$ and a
complex $\cpx{X_2}$ in $\Cb{\pmodcat{A_2}}$ together with an
isomorphism $\cpx{h}: A_0\otimesP_{A_1}\cpx{X_1}\ra
A_0\otimesP_{A_2}\cpx{X_2}$ in $\C{A_0}$, we define a complex
$M(\cpx{X_1}, \cpx{X_2}, \cpx{h}):=(M(X_1^i, X_2^i, h^i), d^i)_{i\in\mathbb{Z}}$, where the differential is induced by the
exact sequence given in Lemma \ref{LemmaPullbackProj}(4). For this complex, we have the following results similar to Lemma \ref{LemmaPullbackProj}.

\begin{Lem} \label{pullbackcpx}Suppose $\cpx{X_1}\in \Cb{\pmodcat{A_1}}$, $\cpx{X_2}\in\Cb{\pmodcat{A_2}}$ and $\cpx{h}: A_0\otimesP_{A_1}\cpx{X_1}\ra
A_0\otimesP_{A_2}\cpx{X_2}$ is an isomorphism in $\C{A_0}$. Then the following hold:

$(1)$ The complex $M(\cpx{X_1}, \cpx{X_2}, \cpx{h})$ is
a bounded complex over $\pmodcat{A}$.

$(2)$ For $i\in\{1,2\}$, there is a natural isomorphism of complexes
$$\cpx{\mu_i}: A_i\otimesP_AM(\cpx{X_1}, \cpx{X_2}, \cpx{h})\lra \cpx{X_i}$$
sending $a_i\otimes (x_1^j, x_2^j)$ to $a_ix_i^j$, and the canonical
projection $\cpx{p_i}: M(\cpx{X_1}, \cpx{X_2}, \cpx{h})\ra \cpx{X_i}$ is equal to
$\lambda_i^*\cpx{\mu_i}$.

$(3)$ There is an exact sequence of complexes of $A$-modules:
$$\xymatrix@M=1.5mm{
0\ar[r] & M(\cpx{X_1}, \cpx{X_2}, \cpx{h}) \ar[r]^(.55){[\cpx{p_1}, \cpx{p_2}]} & \cpx{X_1}\oplus \cpx{X_2} \ar[r]^{\smallvec{\Faf{\pi_1}\cpx{h}}{\Faf{-\pi_2}}} & A_0\otimesP_{A_2}\cpx{X_2}\ar[r] & 0,
}$$
whers $\cpx{p_i}$ is induced by the canonical projection $p_i$ for
$i=1,2$.

$(4)$ Set $\cpx{X}:=M(\cpx{X_1}, \cpx{X_2}, \cpx{h}) $ and $\cpx{X_0}:=A_0\otimesP_{A_2}\cpx{X_2}$.
If $ \Hom_{\K{A_0}}(\cpx{X_0}, \cpx{X_0}[-1])=0$,
then there exists a pullback diagram of algebras:
$$\xymatrix@M=1mm{
\End_{\K{A}}(\cpx{X})\ar[r]^{\epsilon_1}\ar[d]_{\epsilon_2} & \End_{\K{A_1}}(\cpx{X_1})\ar[d]^{\eta_1}\\
\End_{\K{A_2}}(\cpx{X_2})\ar[r]^{\eta_2} & \End_{\K{A_0}}(\cpx{X_0}), \\
}$$
where $\epsilon_1, \epsilon_2, \eta_1$ and $\eta_2$ are homomorphisms of algebras, determined by $\cpx{p_1}$, $\cpx{p_2}$, $\pi_1^*\cpx{h}$ and $\pi_2^*$, respectively.
\end{Lem}

{\it Proof.}
The statements (1)-(3) follow immediately from the definition of
$M(\cpx{X_1}, \cpx{X_2}, \cpx{h})$ and Lemma
\ref{LemmaPullbackProj}(1)-(3). Now, we prove (4).  Since $\cpx{X}\in\Kb{\pmodcat{A}}$, it follows from the triangle
$$\xymatrix@M=1mm{
 \cpx{X} \ar[r]^(.4){[\cpx{p_1}, \cpx{p_2}]} & \cpx{X_1}\oplus \cpx{X_2} \ar[r]^(.55){\smallvec{\Faf{\pi_1}\cpx{h}}{\Faf{-\pi_2}}} & \cpx{X_0}\ar[r] & \cpx{X}[1]
}$$
in $\Db{A}$  that the following long sequence is exact:
$$
\cdots\ra\Hom_{\K{A}}(\cpx{X}, \cpx{X_0}[-1])\ra \Hom_{\K{A}}(\cpx{X}, \cpx{X})\ra\Hom_{\K{A}}(\cpx{X}, \cpx{X_1}\oplus \cpx{X_2})\ra \Hom_{\K{A}}(\cpx{X}, \cpx{X_0})\ra \cdots. $$
Note that $\Hom_{\K{A}}(\cpx{X}, \cpx{X_i}[j])\simeq \Hom_{\K{A_i}}(A_i\otimesP_A\cpx{X},\cpx{X_i}[j])\simeq \Hom_{\K{A_i}}(\cpx{X_i},\cpx{X_i}[j])$ for $j\in \mathbb{Z}$, it follows from the assumption in (4) that $\Hom_{\K{A}}(\cpx{X},\cpx{X_0}[-1])=0$ and the above sequence is then isomorphic to
$$\xymatrix@M=1mm{
0\ar[r] & \End_{\K{A}}(\cpx{X})\ar[r]^(.35){[\epsilon_1, \epsilon_2]} & \End_{\K{A_1}}(\cpx{X_1})\oplus\End_{\K{A_2}}(\cpx{X_2})\ar[r]^(.65){\smallvec{\eta_1}{ -\eta_2}} &\End_{\K{A_0}}(\cpx{X_0}).
}$$
This proves (4).
$\square$

\section{Derived equivalences for Milnor squares of algebras\label{sect3}}

In this section, we first state and prove our main result, Theorem \ref{TheoMain}, on general derived equivalences, and then turn to almost $\nu$-stable derived equivalences (see Corollary \ref{TheoAlmostNuStable}). These derived equivalences induce stable equivalences of Morita type (see \cite{HuXi2010}), while the latter is of interest in an approach to Brou\'e's abelian defect group conjecture (see \cite{HuXi2014-preprint, Rouquier2006}).

\subsection{General result\label{sect3.1}}

The main result of this paper is the following theorem.

\begin{Theo}
Suppose that  $ A_1\lraf{\pi_1} A_0\llaf{\pi_2} A_2$ are homomorphisms of Artin algebras with $\pi_1$ surjective.
Let $\cpx{T_i}$ be a basic, radical tilting complex over $A_i$ with $B_i:=\End_{\Kb{A_i}}(\cpx{T_i})$ for $0\le i\le 2$. If $\cpx{T_0}$ is a direct sum of shifts of
projective $A_0$-modules and there is an isomorphism $A_0\otimesP_{A_i}\cpx{T_i}\simeq
\cpx{T_0}$ of complexes for $i=1,2$, then there exist homomorphisms $ B_1\lraf{\eta_1}B_0\llaf{\eta_2}B_2$ of Artin algebras with $\eta_1$ surjective such that
the pullback algebra $B$ of $\eta_1$ and $\eta_2$ is derived equivalent to the
pullback algebra $A$ of $\pi_1$ and $\pi_2$.
 \label{TheoMain}
\end{Theo}

Thus, it follows immediately from derived invariants that the algebras $A$ and $B$ in Theorem \ref{TheoMain} share many common properties. For instance, they have the same Hochschild (co)homology rings, Coxeter polynomials, and algebraic $K$-theory. For a list of derived invariants, see, for example, \cite{Xi2017} and the references therein.

Remark that if $A_0$ is a product of local algebras, or a self-injective algebra with radical-square zero, then every tilting complex over $A_0$ is  a direct sum of shifts of
projective $A_0$-modules.

To prove Theorem \ref{TheoMain}, we first show the
following lemma.

\begin{Lem}
Let $f: \Lambda\ra \Gamma$ be a surjective homomorphism
between Artin algebras $\Lambda$ and $\Gamma$. If $\cpx{T}$ is a
basic, radical tilting complex over $\Lambda$ such that
 $\Gamma\otimesP_{\Lambda}\cpx{T}$ is a basic tilting complex over $\Gamma$ of the form $\bigoplus_{i=1}^rX_i[n_i]$, where $\{X_1,\cdots, X_r\}$
 is a complete set of non-isomorphic indecomposable projective
 $\Gamma$-modules, then the induced
morphism   $$\Hom_{\K{\Lambda}}(\cpx{T}, f^*):
\Hom_{\K{\Lambda}}(\cpx{T}, \cpx{T})\lra \Hom_{\K{\Lambda}}(\cpx{T},
{}_{\Lambda}\Gamma\otimesP_{\Lambda}\cpx{T})$$ is surjective.
 \label{LemSurjEndoAlg}
\end{Lem}

{\it Proof.}
By Lemma \ref{LemSurjAlgHom}, we can assume that there
is a complete set $\{V_1,\cdots, V_r, U_1,\cdots, U_s\}$ of pairwise
non-isomorphic indecomposable projective $\Lambda$-modules such that
$\Gamma\otimes_{\Lambda}V_i\simeq X_i$ for all $i=1,\cdots, r$, and
that $\Gamma\otimes_{\Lambda}U_i=0$ for all $i=1,\cdots, s$. Set
$U:=\bigoplus_{i=1}^sU_i$. By our assumption, $\nTP{\Gamma\otimes_{\Lambda}\cpx{T}}{X_i}=1$ for all $1\leq i\leq r$.  This implies that $\nTP{\cpx{T}}{V_i}=1$ for $1\leq i\leq r$.  So, by Lemma \ref{LemFormofTiltingComp},  we can write $\cpx{T}$
as
$$\cpx{T}:=\cpx{U}\oplus \cpx{V_1}\oplus\cdots\oplus\cpx{V_r}$$
such that $\cpx{U}\in\Kb{\add(U)}$, and $\nTP{\cpx{V_i}}{V_j}=1$ for
$i=j$ and zero otherwise. Thus
$\Gamma\otimes_{\Lambda}\cpx{U}=0$ and
$\Gamma\otimesP_{\Lambda}\cpx{V_i}\simeq
(\Gamma\otimes_{\Lambda}V_i)[n_i]\simeq X_i[n_i]$ for some integer $n_i$. It is sufficient
to prove that
 $$\Hom_{\K{\Lambda}}(\cpx{T}, f^*): \Hom_{\K{\Lambda}}(\cpx{T},
\cpx{V_i})\lra \Hom_{\K{\Lambda}}(\cpx{T},
{}_{\Lambda}\Gamma\otimesP_{\Lambda}\cpx{V_i})$$ is surjective for
all $i=1,\cdots, r$.

In the following, we set $\Sigma:=\End_{\K{\Lambda}}(\cpx{T})$. Since
$$\begin{array}{rl}
\Hom_{\Sigma}\big(\Hom_{\K{\Lambda}}(\cpx{T},\cpx{U}),
\Hom_{\K{\Lambda}}(\cpx{T}, \Gamma\otimesP_{\Lambda}\cpx{V_i})\big)& \simeq
\Hom_{\K{\Lambda}}\big(\cpx{U}, \Gamma\otimesP_{\Lambda}\cpx{V_i})\big)\\ & \simeq
\Hom_{\K{\Gamma}}(\Gamma\otimesP_{\Lambda}\cpx{U},\Gamma\otimesP_{\Lambda}\cpx{V_i})=0,\end{array}$$
the $\Sigma$-module $\Hom_{\K{\Lambda}}(\cpx{T}, \Gamma\otimesP_{\Lambda}\cpx{V_i})$ has no
composition factors in $\add\!\big(\tp(\Hom_{\K{\Lambda}}(\cpx{T},
\cpx{U}))\big)$.

Now, for $1\leq k\leq r$,  let $S_k$ denote the top of $V_k$, and  $\bar{S}_k$ the top of the $\Sigma$-module
$\Hom_{\K{\Lambda}}(\cpx{T}, \cpx{V_k})$. Let $G:
\Db{\Lambda}\ra\Db{\Sigma}$ be  the derived equivalence induced by
$\cpx{T}$. Then, by the proof of Lemma \ref{LemsimpleTosimple}, we
have $G(S_k)\simeq \bar{S}_k[-n_k]$ for $1\le k\le r$. Since
$\Gamma\otimesP_{\Lambda}\cpx{T}$ is
a tilting complex over $\Gamma$,   $$\Hom_{\K{\Lambda}}(\cpx{T}, (\Gamma\otimesP_{\Lambda}\cpx{V_k})[n])\simeq \Hom_{\K{\Gamma}}(\Gamma\otimesP_{\Lambda}\cpx{T}, (\Gamma\otimesP_{\Lambda}\cpx{V_k})[n])=0$$
for all $1\leq k\leq r$ and  all $n\neq 0$, and
$\Hom_{\K{\Lambda}}(\cpx{T}, \Gamma\otimesP_{\Lambda}\cpx{V_k})\simeq G(\Gamma\otimesP_{\Lambda}\cpx{V_k}) $ for all
$1\leq k\leq r$. Hence
$$\begin{array}{rl}
\Hom_{\Sigma}(\Hom_{\K{\Lambda}}(\cpx{T}, \Gamma\otimesP_{\Lambda}\cpx{V_i}), \bar{S}_k) & \simeq \Hom_{\Db{\Sigma}}(G(\Gamma\otimesP_{\Lambda}\cpx{V_i}), G(S_k[n_k]))\\
\noalign{\smallskip}
  & \simeq \Hom_{\Db{\Lambda}}(\Gamma\otimesP_{\Lambda}\cpx{V_i},  S_k[n_k])\\
  & \simeq \Hom_{\Db{\Lambda}}(X_i[n_i], S_k[n_k])
  \end{array}$$
is zero for all $k\neq i$, and is one-dimensional over
$\End_{\Lambda}(S_k)$ for $k=i$. Hence the top of the
$\Sigma$-module $\Hom_{\K{\Lambda}}(\cpx{T}, \Gamma\otimesP_{\Lambda}\cpx{V_i})$ is
$\bar{S}_i$, and  there is
projective cover
$$\epsilon: \Hom_{\K{\Lambda}}(\cpx{T}, \cpx{V_i})\lra \Hom_{\K{\Lambda}}(\cpx{T}, \Gamma\otimesP_{\Lambda}\cpx{V_i}).$$
Clearly, such an epimorphism is given by
$\Hom_{\K{\Lambda}}(\cpx{T},\cpx{g})$ for some morphism $\cpx{g}:
\cpx{V_i}\ra \Gamma\otimesP_{\Lambda}\cpx{V_i}$. By Lemma \ref{LemmaAlghomForComplex}(2), there is
a morphism $\cpx{u}$ from $\Gamma\otimesP_{\Lambda}\cpx{V_i}$ to
$\Gamma\otimesP_{\Lambda}\cpx{V_i}$ such that $\cpx{g}=f^*\cpx{u}$.
It follows that $$\epsilon=\Hom_{\K{\Lambda}}(\cpx{T}, f^*)\cdot
\Hom_{\K{\Lambda}}(\cpx{T}, \cpx{u}).$$ Hence the endomorphism
$\Hom_{\K{\Lambda}}(\cpx{T}, \cpx{u})$ of the $\Sigma$-module
$\Hom_{\K{\Lambda}}(\cpx{T}, \Gamma\otimesP_{\Lambda}\cpx{V_i})$  is
surjective, and therefore is an isomorphism. Consequently, $\Hom_{\K{\Lambda}}(\cpx{T}, f^*)$ is surjective. This
finishes the proof.
$\square$

\smallskip
{\bf Proof of Theorem \ref{TheoMain}.} We have the following pullback
diagram of homomorphisms of algebras:
$$\xymatrix@M=0.3mm{
  A\ar[r]^{\lambda_1}\ar[d]_{\lambda_2} & *+[r]{A_1}\ar[d]^{\pi_1}\\
  *+[r]{A_2} \ar[r]^{\pi_2} & *+[r]{A_0}.
}$$
By the assumptions of Theorem \ref{TheoMain}, the tilting complex $\cpx{T_0}$ is of the form
$\cpx{T_0}=\bigoplus_{i=1}^m U_i[n_i]$
with $U_i$ projective $A_0$-modules such that $n_i\neq n_j$ whenever $i\neq j$. Thus $\Hom(U_i[n_i], U_j[n_j])=0$ for all $i\neq j$, and $\bigoplus_{i=1}^mU_i$ is a basic,
projective generator for $\modcat{A_0}$.

Recall from Subsection \ref{SectionPullback} that $\pmodcat{A_j}=\add(P_j\oplus Q_j)$ for $j=1,2$, where $A_0\otimes_{A_j}P_j=0$  and $A_0\otimes_{A_j}Y\neq 0$ for each indecomposable direct summand $Y$ of $Q_j$.  Let $\{V_1, \cdots, V_r\}$ and $\{W_1, \cdots, W_s\}$ be  complete sets of pairwise
non-isomorphic indecomposable projective modules in $\add(Q_1)$ and $\add(Q_2)$, respectively. Since $\cpx{h_i}: A_0\otimesP_{A_i}\cpx{T_i}\simeq \cpx{T_0}$ in $\C{A_0}$ for $i=1, 2$, and since each indecomposable projective  $A_0$-module occurs in $\cpx{T_0}$ only once,  we deduce $\nTP{\cpx{T_1}}{V_i}=1=[\cpx{T_2}: W_j]$ for all $i, j$.
By Lemma \ref{LemFormofTiltingComp}, we can write
$$\cpx{T_1}=\cpx{P_1}\oplus \cpx{V_1}\oplus\cdots\oplus\cpx{V_r}\quad\mbox{ and }\quad \cpx{T_2}= \cpx{P_2}\oplus \cpx{{W_1}}\oplus\cdots\oplus \cpx{{W_s}}, $$
such that

(1) $\cpx{P_i}\in\Kb{\add(P_i)}$, and $\add(\cpx{P_i})$ generates
$\Kb{\add(P_i)}$ as a triangulated category for $i=1, 2$, and

(2) $\nTP{\cpx{V_i}}{V_j}=\delta_{ij}$ and $\nTP{\cpx{W_k}}{W_l}=\delta_{kl}$, where $\delta_{ij}$ is the Kronecker symbol.

Note that $A_0\otimes_{{A_1}}P_1=0$ and
$A_0\otimesP_{{A_1}}\cpx{V_i}=(A_0\otimes_{{A_1}}V_i)[n_{{V_i}}]$
for some integer $n_{{V_i}}$ with $1\le i\le r$. By assumption, we have an isomorphism of complexes:
 $$\cpx{h_1}:\quad \bigoplus_{i=1}^r(A_0\otimes_{{A_1}}V_i)[n_{{V_i}}]\simeq\bigoplus_{i=1}^mU_i[n_i].$$
This gives rise to a partition $\sigma=\{\sigma_1,\cdots,\sigma_m\}$ of $\{1,\cdots, r\}$ with $\sigma_i:=\{j\mid n_{{V_j}}=n_i\}$. Now we define
$$V_{\sigma_i}:=\bigoplus_{j\in\sigma_{i}}V_j, \quad\mbox{ and }\quad \cpx{V_{\sigma_i}}:=\bigoplus_{j\in\sigma_{i}}\cpx{V_j}.$$
for all $1\le i\le m$. This partition means that we collect terms of the left-hand side of $\cpx{h_1}$ according to the position $n_i$ of terms in $\cpx{T_0}$. Thus, $$A_0\otimesP_{{A_1}}\cpx{T_1}=\big(A_0\otimesP_{A_1}\cpx{P_1}\big)\oplus \big(\bigoplus_{i=1}^r A_0\otimesP_{{A_1}}\cpx{V_i}\big)=\bigoplus_{i=1}^mA_0\otimesP_{{A_1}}\cpx{V_{\sigma_i}}.$$
Since $\Hom(U_i[n_i], U_j[n_j])=0$ for all $i\neq j$, the
isomorphism $\cpx{h_1}: A_0\otimesP_{{A_1}}\cpx{T_1}\ra \cpx{T_0}$
can be rewritten as
$${\rm diag}[g_1, \cdots, g_m]:\,\,\bigoplus_{i=1}^m(A_0\otimesP_{{A_1}}\cpx{V_{\sigma_i}})\lra \bigoplus_{i=1}^mU_i[n_i]=\cpx{T_0},$$ where $g_i: A_0\otimesP_{{A_1}}\cpx{V_{\sigma_i}}\ra
U_i[n_i]$ is an isomorphism in $\C{A_0}$ for all $i$.
By repeating the above procedure, we get a partition $\tau:=\{\tau_1,\cdots,\tau_m\}$ of $\{1,\cdots,s\}$ with
$$\tau_i:=\{k\in\{1, \cdots, s\}\mid A_0\otimesP_{A_2}\cpx{W_k}\simeq (A_0\otimes_{A_2}W_k)[n_{W_k}]\mbox{ and }n_{W_k}=n_i\}.$$
Define
$$W_{\tau_i}:=\bigoplus_{k\in\tau_{i}}W_k, \quad\mbox{ and }\quad \cpx{W_{\tau_i}}:=\bigoplus_{k\in\tau_{i}}\cpx{W_k}.$$
The isomorphism $\cpx{h_2}$ can be rewritten as
  $${\rm diag}[f_1, \cdots, f_m]:\,\, \bigoplus_{i=1}^mA_0\otimesP_{{A_2}}\cpx{W_{\tau_i}} \lra \bigoplus_{i=1}^mU_i[n_i]=\cpx{T_0}.$$
Now, we define
$\cpx{T}:=M(\cpx{T_1},\cpx{T_2}, \cpx{h_1}\cpx{h_2}^{-1})$, that is, $\cpx{T}=M(\cpx{P_1}, 0, 0)\oplus M(0, \cpx{P_2}, 0)\oplus \bigoplus_{i=1}^m M(\cpx{V_{\sigma_i}}, \cpx{W_{\tau_i}}, g_if_i^{-1}).$ In the sequel, we shall show that $\cpx{T}$ is a tilting complex over $A$.

First, we show that $\add(\cpx{T})$ generates $\Kb{\pmodcat{A}}$ as a triangulated
category.

For simplicity, we write $\cpx{Z_i}$ for $M(\cpx{V_{\sigma_i}}, \cpx{W_{\tau_i}}, g_if_i^{-1})$ for $1\le i\le m$.   By definition, for each integer $k$, $Z_i^k := M(V_{\sigma_i}^k, W_{\tau_i}^k, g_if_i^{-1})$. For $k\neq -n_i$,
the term $V_{\sigma_i}^k$ is in $\add(P_1)$, and the term $W_{\tau_i}^k$ is in
$\add(P_2)$. Hence
$A_0\otimes_{A_1}V_{\sigma_i}^k=0=A_0\otimes_{A_2}W_{\tau_i}^k$, and $Z_i^k\simeq {}_AV_{\sigma_i}^k\oplus
{}_AW_{\tau_i}^k\in\add(F_1\oplus F_2)$ for all $k\neq -n_i$. Since
$V_{\sigma_i}$ is a direct summand of $V_{\sigma_i}^{-n_i}$ and since $W_{\tau_i}$ is a direct
summand of $W_{\tau_i}^{-n_i}$, the $A$-module $M(V_{\sigma_i}, W_{\tau_i},
g_if_i^{-1})$ is a direct summand of $Z_i^{-n_i}$.
By Lemma \ref{LemAddP}(1), the functor $_A(-): \add(P_1)\ra
\add(F_1)$ is an equivalence, and consequently induces a triangle
equivalence between $\Kb{\add(P_1)}$ and $\Kb{\add(F_1)}$. Since
$\add(\cpx{P_1})$ generates $\Kb{\add(P_1)}$ as a triangulated
category, $\add(M(\cpx{P_1}, 0, 0))=\add({}_A\cpx{P_1})$
generates $\Kb{\add(F_1)}$ as a triangulated category. Similarly,
$\add(M(0,\cpx{P_2}, 0))$ generates $\Kb{\add(F_2)}$ as a
triangulated category. As all terms $Z_i^k$ with $k\neq -n_i$
are in $\add(F_1\oplus F_2)$, the term $Z_i^{-n_i}$ is in the
triangulated full subcategory of $\Kb{\pmodcat{A}}$ generated by
$\add(\cpx{T})$. Thus, the module $F_1\oplus F_2\oplus
(\bigoplus_{i=1}^mZ_i^{-n_i})$ is in the triangulated full
subcategory of $\Kb{\pmodcat{A}}$ generated by $\add(\cpx{T})$.  By Lemma \ref{LemAddQ}(2), the direct sum
$$\bigoplus_{i=1}^mM(V_{\sigma_i}, W_{\tau_i}, g_if_i^{-1})$$
is a basic, additive generator of $F_3$. Recall that   $M(V_{\sigma_i}, W_{\tau_i}, g_if_i^{-1})$ is a direct
summand of $Z_i^{-n_i}$ for all $1\le i\le m$.
It follows that $F_1\oplus F_2\oplus F_3$ is generated by $\add(\cpx{T})$ in the triangulated full subcategory of $\Kb{\pmodcat{A}}$ generated by
$\add(\cpx{T})$.  As $F_1\oplus F_2\oplus F_3$ is an additive generator of $\pmodcat{A}$, $\add(\cpx{T})$ generates
$\Kb{\pmodcat{A}}$ as a triangulated category.

Next, we prove that $\cpx{T}$ is self-orthogonal, that is, $\Hom_{\K{A}}(\cpx{T},\cpx{T}[n])=0$ for all $n\ne 0$.

By the construction
of $\cpx{T}$, there is an exact sequence of complexes of $A$-modules (see Lemma \ref{pullbackcpx}(3)):
$$\xymatrix@M=1.5mm{
0\ar[r] & \cpx{T} \ar[r]^(.40){[\cpx{p_1}, \cpx{p_2}]} & \cpx{T_1}\oplus \cpx{T_2} \ar[r]^(0.60){\smallvec{\Faf{\pi_1}\cpx{h_1}}{\Faf{-\pi_2}\cpx{h_2}}} & T_0\ar[r] & 0,
}$$
which yields a triangle in $\Db{A}$.
Applying $\Hom_{\Db{A}}(\cpx{T}, -)$ to this triangle, we get the following commutative diagram with exact rows for each integer $i$:
$$\xymatrix@M=1mm@C=5mm{
 & \Hom_{\Db{A}}(\cpx{T}, \cpx{T_0}[i-1]) \ar[r]\ar[d]^{\simeq}  & \Hom_{\Db{A}}
 (\cpx{T}, \cpx{T}[i])\ar[r]\ar[d]^{\simeq} & \Hom_{\Db{A}}(\cpx{T}, \bigoplus_{k=1}^2\cpx{T_k}[i])\ar[d]^{\simeq}\\
 & \Hom_{\K{A}}(\cpx{T}, \cpx{T_0}[i-1]) \ar[r]\ar[d]^{\simeq}  & \Hom_{\K{A}}
 (\cpx{T}, \cpx{T}[i])\ar[r]\ar@{=}[d] & \bigoplus_{k=1}^2\Hom_{\K{A}}(\cpx{T}, \cpx{T_k}[i])\ar[d]^{\simeq}\\
(**)\quad& \Hom_{\K{A_0}}(\cpx{T_0}, \cpx{T_0}[i-1]) \ar[r]  &
\Hom_{\K{A}}(\cpx{T}, \cpx{T} [i])\ar[r]&
\bigoplus_{k=1}^2\Hom_{\K{A_k}}
(\cpx{T_k}, \cpx{T_k}[i])\\
}$$ Here we use the following natural isomorphisms: $$\Hom_{\K{A}}(\cpx{T}, \cpx{T_k}[i])\simeq \Hom_{\K{A_k}}(A_k\otimes_A\cpx{T}, \cpx{T_k}[i])\simeq \Hom_{\K{A_k}}(\cpx{T_k}, \cpx{T_k}[i])$$ for $0\le k\le 2,$ where the last isomorphism is due to Lemma \ref{pullbackcpx}(2). Since $\Hom_{\K{A_0}}(\cpx{T_0}, \cpx{T_0}[i-1])=0$ for
all $i\neq 1$ and since $\Hom_{\K{A_k}}
(\cpx{T_k}, \cpx{T_k}[i])=0$ for all
$i\neq 0$ and all $0\le k\le 2$, we have $\Hom_{\K{A}}(\cpx{T},\cpx{T}[i])=0$ for all
$i\neq 0,1$. It follows from Lemma \ref{LemSurjEndoAlg}
 that the morphism $\eta_1: \Hom_{\K{A_1}}(\cpx{T_1},
\cpx{T_1}) \ra \Hom_{\K{A_0}}(\cpx{T_0}, \cpx{T_0})$ determined by
$\pi_1^*\cpx{h_1}$ is surjective. Consequently, from the long exact sequence
$(**)$, we get $\Hom_{\K{A}}(\cpx{T},\cpx{T}[1])$ = $0$. Thus $\cpx{T}$ is self-orthogonal.
Altogether, we have shown that $\cpx{T}$ is a
tilting complex over $A$.

To finish the proof of Theorem \ref{TheoMain}, we consider the endomorphism algebra of $\cpx{T}$. By Lemma \ref{pullbackcpx}(4), there exists a pullback diagram of homomorphisms of algebras:
 $$\xymatrix@M=0.5mm{
  \End_{\K{A}}(\cpx{T})\ar[r]^{\epsilon_1}\ar[d]^{\epsilon_2} & \End_{\K{A_1}}(\cpx{T_1})\ar@{->>}[d]^{\eta_1}\\
  \End_{\K{A_2}}(\cpx{T_2}) \ar[r]^{\eta_2} &
  \End_{\K{A_0}}(\cpx{T_0}),
}$$
where $\eta_1$ and $\eta_2$ are
determined by $\pi_1^*\cpx{h_1}$ and
$\pi_2^*\cpx{h_2}$, respectively, and where $\epsilon_1$ and $\epsilon_2$ are
 determined by the projections from $\cpx{T}$ to $\cpx{T_1}$ and $\cpx{T_2}$, respectively.
This completes the proof of Theorem \ref{TheoMain}.
$\square$

An immediate consequence of Theorem \ref{TheoMain} is the following result.

\begin{Koro} Let $A$ be an Artin algebra and $T$ a basic, radical tilting complex over $A$. Suppose that $I$ is an ideal in $A$ such that $\rad(A)\subseteq I$, $\Hom_{\Kb{A}}\big(\cpx{T},I\cpx{T}[i]\big)$ = $0$ for all $i\ne 0$ and $\Hom_{\Kb{A}}(\cpx{T}/I\cpx{T}, (\cpx{T}/I\cpx{T})[-1])$ = $0$. Let $B:=\End_{\Kb{A}}(\cpx{T})$ and $J$ be the ideal of $B$ consisting of all those endomorphisms of $\cpx{T}$ that factorize through the injection $I\cpx{T}\ra \cpx{T}$. If $\cpx{T}/I\cpx{T}$ is a basic, radical complex, then
the algebras $$ \Lambda:=\{(a,a')\in A\times A\mid a-a'\in I\} \; \mbox{  and  } \; \Gamma:=\{(b,b')\in B\times B\mid b-b'\in J\}$$ are derived equivalent. \label{corollary1}
\end{Koro}

{\it Proof.} By the assumptions on $I$ and \cite[Theorem 4.2]{HuXi2013}, we see that the complex $\cpx{T}/I\cpx{T}$ is a tilting complex over $A/I$ and induces a derived equivalence between $A/I$ and $B/J$. Since the algebra $A/I$ is semisimple, the complex $\cpx{T}/I\cpx{T}$ satisfies the conditions of Theorem \ref{TheoMain} for $\cpx{T_0}$. Thus the pullback algebras of $A\ra A/I\la A$ and $B\ra B/J\la B$ are derived equivalent, that is, $\Lambda$ and $\Gamma$ are derived equivalent. $\square$

\subsection{Special case: iterated almost $\nu$-stable derived equivalences\label{sect3.2}}

A special class of derived equivalences is the one of almost $\nu$-stable derived equivalences which induce stable equivalences of Morita type, while such stable equivalences play a significant role in an approach to Brou\'e's abelian defect group conjecture (see \cite{Rouquier2006, HuXi2014-preprint}). Thus it is quite natural to ask if almost $\nu$-stable derived equivalences can be constructed from Milnor squares. In this section, we show that it is the case for finite-dimensional algebras over an algebraically closed field (see Corollary \ref{TheoAlmostNuStable}).

Throughout this section all algebras are finite-dimensional over a fixed field.

Let $F:\Db{A}\ra\Db{B}$ be a derived equivalence of algebras $A$ and $B$. Suppose that $\cpx{Q}$ and $\cpx{\bar{Q}}$ are radical tilting complexes associated to $F$ and the quasi-inverse $F^{-1}$ of $F$, respectively. By applying the shift function if necessary, we may assume that $\cpx{Q}$ is of the form
$$0\lra Q^{-n}\lra\cdots \lra Q^{-1}\lra Q^0\lra 0$$
and $\cpx{\bar{Q}}$ is of the form
$$0\lra {\bar Q}^0\lra {\bar Q}^1\lra\cdots\lra {\bar Q}^n\lra 0.$$
Let $Q:=\bigoplus_{i=1}^nQ^{-i}$ and $\bar{Q}:=\bigoplus_{i=1}^n\bar{Q}^n$. The derived equivalence $F$ is called {\em almost $\nu$-stable} provided that $\add({}_AQ)=\add(\nu_AQ)$ and $\add({}_B\bar{Q})=\add(\nu_B\bar{Q})$. The composite of finitely many almost $\nu$-stable derived equivalences or their quasi-inverses is called an {\em iterated almost $\nu$-stable derived equivalence}. Such a derived equivalence of finite-dimensional algebras over a field always induces a stable equivalence of Morita type (see \cite{HuXi2010} and \cite{Hu2012}).

A module $P\in \modcat{A}$ is said to be $\nu$-\emph{stably projective} if $\nu_A^iP$ is projective for all $i\geq 0$, where $\nu_A$ is the Nakayama functor $D\Hom_A(-,A)\simeq D(A)\otimes_A-: \modcat{A}\ra \modcat{A}$. We denote by $\stp{A}$ the full subcategory of $\pmodcat{A}$ consisting of all $\nu$-stably projective $A$-modules.

For finite-dimensional algebras, Theorem \ref{TheoMain} can be strengthened as the following corollary which is the main result in this subsection.

\begin{Koro}\label{TheoAlmostNuStable}
Keep the assumptions in Theorem \ref{TheoMain}, and further assume the following conditions:

$(1)$ $A_0, A_1$ and $A_2$ are finite-dimensional algebras over an algebraically closed field $k$.

$(2)$ $\cpx{T_1}$ and $\cpx{T_2}$ induce iterated almost $\nu$-stable derived equivalences.

$(3)$ $\cpx{T_0}$ is a stalk complex concentrated in degree zero.

\noindent Then the derived equivalence between the pullback algebras in Theorem \ref{TheoMain} is  iterated almost $\nu$-stable.
 \end{Koro}

Thus the pullback algebras in Corollary \ref{TheoAlmostNuStable} have many common nice properties: the same global, finitistic and dominant dimensions, and the same numbers of non-isomorphic, non-projective simple modules, that is the Auslander-Reiten conjecture holds true for the two stably equivalent algebras (see \cite{HuXi2010, HuXi2014-preprint}).

For the proof of Corollary \ref{TheoAlmostNuStable}, we have to prepare a few lemmas. Recall that $S_X$ denotes the top of an indecomposable projective module $X$.

\begin{Lem} If $P$ is an indecomposable module in $\stp{A}$, then there is an exact sequence of $A$-modules
$$(\star)\quad 0\lra R_P\lra \nu_{A}S_P\lra S_{\nu P}\lra 0$$
such that the composition factors of $R_P$ are of the form $S_X$ for some indecomposable projective $X\notin\stp{A}$. \label{compositionfactor}
\end{Lem}

{\it Proof.} Since $S_P$ is the top of $P$, the module $\nu_{A}S_P$ is a quotient of $\nu_{A}P\in \stp{A}$, while $\nu_{A}P$ is an indecomposable projective module in $\stp{A}$ and has $S_{\nu_AP}$ as its top. Thus $\nu_{A}S_P$ is an indecomposable module with a simple top $S_{\nu_AP}$.
Hence there is an exact sequence of $A$-modules:
$$ 0\lra R_P\lra \nu_{A}S_P\lra S_{\nu P}\lra 0.$$

For each indecomposable module $Y\in\stp{A}$, the multiplicity of $S_Y$ as a composition factor of $\nu_{A}S_P$ is the length of $\Hom_{A}(Y, \nu_{A}S_P)$ as an $\End_A(S_Y)$-module. However,
$$\Hom_A(Y, \nu_AS_P)\simeq\Hom_A(Y, D(A)\otimes_{A}S_P)\simeq \Hom_{A}(Y, D(A))\otimes_{A}S_P\simeq D(Y)\otimes_{A}S_P\simeq\Hom_{A}(\nu_{A}^{-1}Y, S_P)$$
is zero if $Y\not\simeq \nu_{A}P$, and has length $1$ if $Y\simeq \nu_{A}P$. Hence $\nu_{A}S_P$ has the composition factor $S_{\nu P}$ at top with $[\nu_{A}S_P:S_{\nu P}]=1$, and  other composition factors of the form $S_X$ with $X$ an indecomposable projective module not in $\stp{A}$.
$\square$

\begin{Lem}{\rm \cite[Theorem 1.1]{Hu2012}} Let $F:\Db{A}\ra\Db{B}$ be a derived equivalence between
algebras $A$ and $B$ over an algebraically closed field, and let $\cpx{T}$ and $\cpx{\bar{T}}$ be tilting complexes
associated to $F$ and $F^{-1}$, respectively. Set $T^{\pm}:=\oplus_{i\ne 0} T^i$ and $\bar{T}^{\pm}:=\oplus_{j\ne 0}{\bar{T}}^j$.
Then the following are equivalent:

$(1)$ The functor $F$ is an iterated almost $\nu$-stable derived equivalence;

$(2)$ $\add(T^{\pm})=\add(\nu_A T^{\pm})$ and $\add(\bar{T}^{\pm})=\add(\nu_B{\bar T}^{\pm})$;

$(3)$ $T^{\pm}\in \stp{A}$ and $\bar{T}^{\pm}\in\stp{B}$;

$(4)$ For each indecomposable projective $A$-module $X\notin \stp{A}$, $F(S_X)$ is
isomorphic in $\Db{B}$ to a simple $B$-module;

$(5)$ For each indecomposable projective $A$-module $X\notin\stp{A}$, there hold
$X\notin \add(T^{\pm})$ and $[U^0:X]=1$, where
$\cpx{U}=(U^i,d_U)$ is the direct sum of all non-isomorphic indecomposable direct summands
of $\cpx{T}$.

Moreover, if one of (1)-(5) is satisfied, then $A$ and
$B$ are stably equivalent of Morita type.\label{huiteratedstaeui}
\end{Lem}

Thus a derived equivalence $F$ is iterated almost $\nu$-stable if and only if so is its quasi-inverse $F^{-1}$ by (2). For the definition of stable equivalences of Morita type, the reader is referred to, for instance, \cite{HuXi2010}.

\begin{Lem}\label{LemmaSimpleUnderAlmosNuStable}
Let $\Lambda$ and $\Gamma$ be algebras over an algebraically closed field and  $F:\Db{\Lambda}\ra\Db{\Gamma}$ be an iterated almost $\nu$-stable derived equivalence. Suppose that $P$ is an indecomposable projective $\Lambda$-module in $\stp{\Lambda}$.

$(1)$ If $F(S_P)$ is isomorphic to a simple $\Gamma$-module $S_{P'}$, then so is $F(S_{\nu_{\Lambda}P})$. Moreover, $P'$ must be in $\stp{\Gamma}$.

$(2)$ If $F(S_P)$ is not isomorphic to a simple $\Gamma$-module, then neither is $F(S_{\nu_{\Lambda}P})$.
\end{Lem}

{\it Proof.}
(1) We may assume that the given derived equivalence $F$ is almost $\nu$-stable with $\cpx{Q}$ and $\cpx{\bar{Q}}$ being radical tilting complexes associated to $F$ and $F^{-1}$, respectively. Let $Q:=\bigoplus_{i>0}Q^{-i}$ and $\bar{Q}:=\bigoplus_{i>0}\bar{Q}^i$. Then, by definition, $\add(\nu_{\Lambda}Q)=\add(Q)$ and $\add(\nu_{\Gamma}\bar{Q})=\add(\bar{Q})$.

By \cite[Lemma 5.2]{HuXi2010}, there is a radical, two-sided tilting complex ${}_{\Gamma}\cpx{\Delta_{\Lambda}}$:
$$0\lra \Delta^0\lra \Delta^1\lra\cdots\lra \Delta^n\lra 0$$
such that $F(\cpx{X})\simeq \cpx{\Delta}\otimesP_{\Lambda}\cpx{X}$ with $\Delta^i\in\add(\bar{Q}\otimes_kQ^*)$ for all $i>0$. Here, $Q^*=\Hom_{\Lambda}(Q,\Lambda)$ is the $\Lambda$-duality of $_{\Lambda}Q$. Let $\cpx{\Theta}:=\HomP_{\Gamma}(\cpx{\Delta},\Gamma)$, the an inverse of $\cpx{\Delta}$. Then the bimodules $\Delta^0$ and $\Theta^0$
define a stable equivalence of Morita type between $\Lambda$ and $\Gamma$ (see the proof of Theorem 5.3 in \cite{HuXi2010}). Here, we stress that $\Delta^0\otimes_{\Lambda}-$ is both a left and right adjoint to $\Theta^0\otimes_{\Gamma}-$. Indeed, $\Theta^0: =\Hom_{\Gamma}(\Delta^0,\Gamma)$ implies that  $\Delta^0\otimes_{\Lambda}-$ is a left adjoint to  $\Theta^0\otimes_{\Gamma}-$. Note that there is an isomorphism $\cpx{\Delta}\simeq\HomP_{\Lambda}(\cpx{\Theta},\Lambda)$ in $\Db{\Gamma\otimes_k\Lambda\opp}$, due to the fact that $\cpx{\Delta}$ is an inverse of $\cpx{\Theta}$. The isomorphism can be regarded as in $\Kb{\Gamma\otimes_k\Lambda\opp}$ by \cite[Lemma 2.1]{HuXi2010}. Since both complexes $\cpx{\Delta}$ and $\HomP_{\Lambda}(\cpx{\Theta},\Lambda)$ are  radical, they are even isomorphic in $\Cb{\Gamma\otimes_k\Lambda\opp}$ by \cite[(b), p.113]{HuXi2010}. It follows that $\Delta^0\simeq\Hom_{\Lambda}(\Theta^0, \Lambda)$ and  $\Delta^0\otimes_{\Lambda}-$ is a  right adjoint to  $\Theta^0\otimes_{\Gamma}-$.

Suppose  $F(S_P)\simeq S_{P'}$ in $\Db{\Gamma}$ for an indecomposable projective $\Gamma$-module $P'$. Then $P'\in \stp{\Gamma}$. In fact, if $P'\notin\stp{\Gamma}$, then
$\Hom_{\Lambda}(P, S_P)\simeq\Hom_{\Db{\Gamma}}(F(P), S_{P'})$
would vanish since $F(P)$ is isomorphic to a complex in $\Kb{\stp{\Gamma}}$ by \cite[Lemma 3.9]{HuXi2010}. This is a contradiction.

To prove (1), we show  $F(S_{\nu P})\simeq S_{\nu P'}$.

Indeed, since $F(S_P)$ is simple, $\Hom_{\Db{\Lambda}}(\cpx{T}, S_P[i])\simeq \Hom_{\Db{\Gamma}}(\Gamma, F(S_P)[i])=0$ for all $i\neq 0$. It follows that
$Q^*\otimes_{\Lambda}S_P\simeq\Hom_{\Lambda}(Q, S_P)=0$, and thus $\Delta^i\otimes_{\Lambda}S_P=0$ for all $i> 0$. Hence
$F(S_P)\simeq\cpx{\Delta}\otimesP_{\Lambda}S_P=\Delta^0\otimes_{\Lambda}S_P\simeq S_{P'}$.

For $P\in \stp{\Lambda}$, there is the following exact sequence of $\Lambda$-modules by Lemma \ref{compositionfactor}:
$$(\star)\quad 0\lra R_P\lra \nu_{\Lambda}S_P\lra S_{\nu P}\lra 0$$
Now, applying $\Delta^0\otimes_{\Lambda}-$ to ($\star$), we get an exact sequence of $\Gamma$-modules
$$(\star\star)\quad 0\lra \Delta^0\otimes_{\Lambda}R_P\lra\Delta^0\otimes_{\Lambda}\nu_{\Lambda}S_P\lra\Delta^0\otimes_{\Lambda}S_{\nu P}\lra 0.$$

Note that $\Delta^0\otimes_{\Lambda}\nu_{\Lambda}S_P\simeq\nu_{\Gamma} (\Delta^0\otimes_{\Lambda}S_{P})$ by a property of stable equivalences of Morita type (see (b) in the proof of \cite[Lemma 3.1]{HuXi2014-preprint}. Note that (b) holds without any additional
assumptions in \cite[Lemma 3.1]{HuXi2014-preprint} because $\Delta^0\otimes_{\Lambda}-$ is both a left and right adjoint to $\Theta^0\otimes_{\Gamma}-$).
Recall that $F(S_P)\simeq \Delta^0\otimes_{\Lambda}S_P\simeq S_{P'}$. It follows that $\nu_{\Gamma} (\Delta^0\otimes_{\Lambda}S_{P})\simeq \nu_{\Gamma}(F(S_P))\simeq \nu_{\Gamma}(S_{P'}).$
Hence $\Delta^0\otimes_{\Lambda}\nu_{\Lambda}S_P\simeq \nu_{\Gamma}(S_{P'}).$ Due to $\Hom_{\Lambda}(Q,S_P)=0$, we get $P\notin\add(Q)$ and $\nu_{\Lambda}P\notin\add(\nu_{\Lambda} Q)=\add(Q)$. This implies $\Hom_{\Lambda}(Q, S_{\nu P})=0$. Hence $\Delta^i\otimes_{\Lambda}S_{\nu P}=0$ for $i\neq 0$ and $F(S_{\nu P})\simeq\cpx{\Delta}\otimes_{\Lambda}S_{\nu P}\simeq\Delta^0\otimes_{\Lambda}S_{\nu P}$. Thus we assume $F(S_{\nu P})=\Delta^0\otimes_{\Lambda}S_{\nu P}\in \modcat{\Gamma}$ and rewrite ($\star\star$) as

$$0\lra \Delta^0\otimes_{\Lambda}R_P\lra \nu_{\Gamma}S_{P'}\lra F(S_{\nu P})\lra 0 $$

Note that both $\nu_{\Gamma}S_{P'}$ and $F(S_{\nu P})$ have a simple top isomorphic to $S_{\nu P'}$ and that $\nu_{\Gamma}S_{P'}$ has other composition factors of the form $S_{X'}$ with $X'\notin\stp{\Gamma}$ indecomposable by Lemma \ref{compositionfactor}. So, to prove that $F(S_{\nu P})$ is simple, we only have to show that $F(S_{\nu P})$ does not have any submodule isomorphic to $S_{X'}$ for all indecomposable projective $\Gamma$-modules $X'\notin \stp{\Gamma}$. This is equivalent to showing
$\Hom_{\Gamma}(S_{X'}, F(S_{\nu P}))=0$ for all indecomposable projective modules $X'\notin\stp{\Gamma}$. Indeed, by definition, $F$ is iterated almost $\nu$-stable if and only if $F^{-1}$ is iterated almost $\nu$-stable. Hence, by Lemma \ref{huiteratedstaeui}(4),
for each indecomposable projective $\Gamma$-module $X'\notin\stp{\Gamma}$, there is an indecomposable projective $\Lambda$-module  $X\notin\stp{\Lambda}$ such that $F(S_X)\simeq S_{X'}$.
Thus $\Hom_{\Gamma}(S_{X'}, F(S_{\nu P}))\simeq \Hom_{\Lambda}(S_X, S_{\nu P})=0$. Consequently, $F(S_{\nu P})$ has a unique composition factor $S_{\nu P'}$, that is, $F(S_{\nu P})\simeq S_{\nu P'}$.

(2) follows from (1).
$\square$

\medskip
\textbf{Proof of Corollary \ref{TheoAlmostNuStable}.} We keep the notations in the proof of Theorem \ref{TheoMain}. The tilting complex $\cpx{T}$ induces a derived equivalence between the pullback algebras.
%By Proposition \ref{huiteratedstaeui},
To prove that $\cpx{T}$ induces an iterated almost $\nu$-stable derived equivalence, we show the following statements:

(a) $T^i\in\stp{A}$ for all $i\neq 0$

In fact, by assumption, the complex $\cpx{T_0}$ is a stalk complex concentrated in degree $0$ and $A_0\otimes_{A_i}\cpx{T_i}\simeq\cpx{T_0}$ for $i=1,2$. It follows that $T_i^m\in\add(P_i)$ for $i=1,2$ and $m\neq 0$, where $P_i$ is as defined in  Subsection \ref{SectionPullback}. Thus, by the construction of $\cpx{T}$, the term $T^m$ is equal to $M(T_1^m, 0, 0)\oplus M(0, T_2^m, 0)$ for $m\neq 0$. By Lemma \ref{huiteratedstaeui}, for $i\in\{1,2\}$, the $A_i$-module $T_i^{\pm}:=\bigoplus_{m\neq 0}T_i^m$ satisfies $\add(\nu_{A_i}T_i^{\pm})=\add(T_i^{\pm})$. It follows from Lemma \ref{LemAddP}(3) that $T^{\pm}:=\bigoplus_{m\neq 0}T^m$ satisfies $\add(\nu_AT^{\pm})=\add(T^{\pm})$. Hence $T^m\in\stp{A}$ for all $m\neq 0$.

(b)
$\nTP{T^0}{X}=1$ for each indecomposable projective $A$-module $X\notin \stp{A}$.

Let $X$ be an indecomposable projective $A$-module and $X\notin \stp{A}$. We need to show  $\nTP{T^0}{X}=1$. Suppose contrarily $\nTP{T^0}{X}=r>1$. Clearly, from the construction of $\cpx{T}$, we have $T^0\simeq M(T_1^0, T_2^0, h^0)$ with  $h^0: A_0\otimes_{A_1}T_1^0\ra A_0\otimes_{A_2}T_2^0$ an isomorphism of $A_0$-modules. Also, $X\simeq M(X_1, X_2, h_X)$ for $X_1=A_1\otimes_AX$, $X_2=A_2\otimes_AX$ and an $A_0$-module isomorphism $h_X: A_0\otimes_{A_1}X_1\ra A_0\otimes_{A_2}X_2$. If $h_X\neq 0$, then $A_0\otimes_{A_1}X_1=A_0\otimes_{A_1}A_1\otimes_AX$ is a direct summand of $A_0\otimes_{A_1}A_1\otimes_AT^0\simeq A_0\otimes_{A_1}T_1^0\simeq T_0^0$ with the multiplicity at least $r$. This contradicts to the assumption that $\cpx{T_0}$ is a basic projective generator of $A_0$-modules. Hence $h_X=0$, $A_0\otimes_{A_1}X_i=0$  for $i=1,2,$ and $X\simeq M(X_1, 0, 0)\oplus M(0, X_2, 0) = X_1\oplus X_2$. It follows that $X_i\in\add(P_i)$ for $i=1,2$, and either $X_1=0$ or $X_2=0$. Without loss of generality, we assume $X_1\neq 0$. Then $[T_1^0: X_1]\geq r$ since $A_1\otimes_AT^0\simeq T_1^0$, and consequently $X_1\in \stp{A_1}$ by Lemma \ref{huiteratedstaeui}(5), and the image of $\tp(X_1)$ of the indecomposable projective $A_1$-module $X_1$ under the derived equivalence induced by $\cpx{T_1}$ is not isomorphic to a simple module by Lemma \ref{LemsimpleTosimple}.

To finish the proof of (b), we show the following:

\smallskip
($\divideontimes$) Let $Y_1$ be the direct sum of all indecomposable projective $A_1$-modules $Y$ such that the image of $\tp(Y)$ under the derived equivalence induced by $\cpx{T_1}$ is not isomorphic to a simple module. Then $Y_1\in\add(P_1)$ and $\add(\nu_{A_1}Y_1)=\add(Y_1)$.

\smallskip
Indeed, let $Y$ be an indecomposable $A_1$-module such that the image of $\tp(Y)$ under the derived equivalence induced by $\cpx{T_1}$ is not isomorphic to a simple module. If $Y$ is a direct summand of $T_1^m$ for some $m\ne 0$, then $Y\in\add(P_1)$ since $T_1^m\in\add(P_1)$. Now, assume that $Y$ only occurs, as a direct summand, in $T_1^0$. Then $[T_1^0: Y]>1$ by Lemma \ref{LemsimpleTosimple}. If $Y\notin\add(P_1)$, then $0\ne A_0\otimes_{A_1}Y$ is a direct summand of $T_0^0$. It follows from $T_0^0\simeq A_0\otimes_{A_1}T_1^0$ with $[T_0^0: A_0\otimes_{A_1}T_1^0]>1$ that $T_0^0=\cpx{T_0}$ is not a basic $A_0$-module. Consequently, $Y$ cannot occur in $T_1^0$, and therefore $Y\in\add(P_1)$ and $Y_1\in \add(P_1)$.
By assumption and Lemma \ref{huiteratedstaeui}(4), $Y_1\in \stp{A_1}$. Now, it follows from Lemma \ref{LemmaSimpleUnderAlmosNuStable}(2) that, for each indecomposable $Y\in\add(Y_1)$, the module $\nu_{A_1}Y$ is again in $\add(Y_1)$. Hence $\add(\nu_{A_1}Y_1)=\add(Y_1)$.

Thus, by ($\divideontimes$) and Lemma \ref{LemAddP}(3), $X=X_1$ lies in $\stp{A}$. This is a contradiction and shows $\nTP{T^0}{X}=1$.

Altogether, we have shown $\add(\nu_AT^{\pm})=\add(T^{\pm})$, $T^{\pm}\in \stp{A}$ and $\nTP{T^0}{X}=1$ for every indecomposable projective $A$-module $X\notin\stp{A}$. Note that if $X\notin\stp{A}$ then $X\notin \add(T^{\pm})$. Now, by Lemma \ref{huiteratedstaeui}(5), $\cpx{T}$ induces an iterated almost $\nu$-stable derived equivalence.
$\square$

\section{Some realizations by quivers with relations \label{sect4}}

In this section, we shall realize the main result, Theorem \ref{TheoMain}, by three ``local" operations on derived equivalent algebras presented by quivers with relations. They are facilitated by gluing vertices, unifying arrows and identifying socle elements. The details are given in
Theorems \ref{TheoremGluing}, \ref{TheoUnifyingArrows} and \ref{theorem-gluesoc}, respectively. Note that these operations can be combined with each other and applied repeatedly.

Let $Q=(Q_0, Q_1)$ be a quiver with $Q_0$ the set of vertices
and $Q_1$ the set of arrows between vertices. For $m>1$, let $Q_m$ be the set of all paths in $Q$ of length $m$. The starting and ending vertices of a path $p$ are denoted by $s(p)$
 and $e(p)$, respectively. As usual, the trivial path corresponding to a vertex
 $i\in Q_0$ is denoted by $e_i$.

We fix a field $k$ and denote by $kQ$ the path algebra of $Q$ over $k$. The composition of two paths  $p$ and $q$ in $kQ$ is written as $pq$
if $e(p)=s(q)$, and zero otherwise. A \emph{ relation} $\omega$
 on $Q$ is a $k$-linear combination of paths: $\omega=\lambda_1p_1+\lambda_2p_2+\cdots+\lambda_np_n$
 with $0\ne\lambda_i\in k$, $e(p_1)=\cdots =e(p_n)$ and
 $s(p_1)=\cdots=s(p_n)$. Here, we assume that the length of each
 $p_i$, that is the number of arrows in $p_i$, is at least 2.  If $n=1$ in $\omega$, then $\omega$ is called a \emph{monomial} relation.

  Let $\rho$ be a set of relations in $kQ$ and $\langle\rho\rangle$ be the ideal
  of $kQ$ generated by $\rho$. Then  an algebra of the form $kQ/\langle\rho\rangle$
  is said to be presented by the quiver $Q$ with relations $\rho$.
  Clearly, $\langle\rho\rangle\subseteq \langle Q_2\rangle$. Note that for any ideal
  $I\subseteq \langle Q_2\rangle$ of $kQ$ such that $kQ/I$ is finite-dimensional, there is a set $\rho$ of relations such that $\langle\rho\rangle=I$.

\subsection{Derived equivalences from gluing vertices}\label{subsection-gluing-vertices}
In this subsection, we shall construct derived equivalences from given ones by gluing
  vertices of quivers. This also gives a way to get derived equivalences for subalgebras from the ones for given algebras.

Let $A=kQ/\langle\rho\rangle$ be a
finite-dimensional algebra over a field $k$. For a subset $X\subseteq Q_0$, we denote by $e_X$ the idempotent element $\sum_{i\in X}e_i$ in $A$.  Let $X$ be a subset of $Q_0$ and $\sigma=\{\sigma_1, \cdots, \sigma_m\}$ be a partition of $X$, that is, $X=\bigcup_i\sigma_i$ and $\sigma_i\cap\sigma_j=\O$ for $i\neq j$.
Let $Q^{\sigma}$ be the quiver obtained from $Q$ by just gluing the vertices in $\sigma_t$
into one vertex, also denoted by $\sigma_t$, for all $t$, and keeping all arrows. Thus the vertex set of $Q^{\sigma}$ is the union of $\{\sigma_1, \sigma_2,\cdots,\sigma_m\}$ and $Q_0\setminus X$, and the arrow set of $Q^{\sigma}$ is $Q_1$. Then there is a natural homomorphism of algebras:
    $$\lambda_{\sigma}: kQ^{\sigma}\lra kQ/\langle\rho\rangle$$
which sends $e_i$ to $e_i$ for $i\not\in X$, $e_{{\sigma}_t}$ to $\sum_{i\in {\sigma_t}}e_i$ for $1\le t\le m$
and preserves all arrows. Clearly, the kernel of $\lambda_{\sigma}$ is contained in
  $\langle Q^{\sigma}_2\rangle$ in $kQ^{\sigma}$.
  Let $\rho^{\sigma}$ be a set of relations on $Q^{\sigma}$ such that $\langle\rho^{\sigma}\rangle=\Ker(\lambda_{\sigma})$.
  The relations $\rho^{\sigma}$ can be obtained in the following way: For each $t$, let $\rho^{\sigma_t}$ be the set of
  relations on $Q^{\sigma}$ consisting of all $\alpha\beta$ with $\alpha$, $\beta$ being arrows such that
  $e(\alpha)$ and $s(\beta)$ are different vertices in $\sigma_t$.
  Then $\rho^{\sigma}=\rho\cup
  \rho^{\sigma_1}\cup\cdots\cup\rho^{\sigma_m}$. The algebra
  $A^{\sigma}:=kQ^{\sigma}/\langle\rho^{\sigma}\rangle$ is called the \emph{$\sigma$-gluing algebra} of
  $A$. The above homomorphism $\lambda_{\sigma}$ induces a homomorphism from $A^{\sigma}$ to $A$, denoted again by $\lambda_{\sigma}$. Observe that $\lambda_{\sigma}: A^{\sigma}\ra A$ is injective, and the
  image of $\lambda_{\sigma}$ is the subalgebra of $A$ generated by all the arrows in
  $Q$, the idempotents
  $e_{{\sigma_1}}, \cdots, e_{{\sigma_m}}$ and $\{e_i\mid i\in Q_0\backslash X\}$.
Note that the Jacobson radicals of $A^{\sigma}$ and $A$ are equal. This construction has been used in the study of
the finitistic dimension conjecture (for example, see \cite{Xi2006}).

Now, we illustrate the above procedure by an example. Let $A$ be a $k$-algebra presented by the quiver $Q$
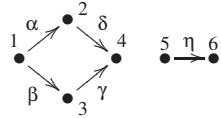
$$\xymatrix@C=4mm@R=3mm{
      & \bullet\ar[rd]^{\delta}^(.1){2}^(.9){4}\\
    \bullet\ar[ru]^{\alpha}^(.1){1}\ar[rd]_{\beta} & & \bullet   & \bullet\ar[r]^{\eta}^(0){5}^(1){6} & \bullet\\
    & \bullet \ar[ru]_{\gamma}_(.1){3} \\
  }$$
with the relation $\alpha\delta-\beta\gamma$. Let
$X:=\{1,2,3,4,5\}$ and $\sigma:=\big\{\{1,2,3\},
\{4,5\}\big\}$ be a partition of $X$. Then the
$\sigma$-gluing algebra $A^{\sigma}$ of $A$ is presented by the quiver
$Q^{\sigma}$
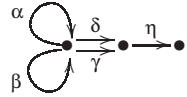
$$\xymatrix@C=5mm@R=8mm{
 {\bullet}\ar@<-2pt>@(l,u)^(.4){\alpha}\ar@<-2pt>@(l,d)_(.4){\beta}\ar@<2pt>[r]^{\delta}\ar@<-2pt>[r]_{\gamma} &\bullet\ar[r]^{\eta} & \bullet\\
}$$
with relations $\{\alpha\delta-\beta\gamma\}\cup \rho^{\sigma_1}\cup\rho^{\sigma_2}=\{\alpha\delta-\beta\gamma,\alpha^2,\alpha\beta,\alpha\gamma, \beta\alpha,\beta^2,\beta\delta,\delta\eta,\gamma\eta\}$.

\medskip
In the following, we shall interpret the procedure of a $\sigma$-gluing as a pullback of algebras. We define
      $$k^X:=\bigoplus_{i\in X}k$$
to be the path algebra of the quiver with isolated
vertices indexed by $X$. Considering $\sigma$ as a set, we
have the algebra $k^{\sigma}$ which is just the $\sigma$-gluing algebra of
$k^{X}$. There is an embedding
$\lambda_{\sigma}: k^{\sigma}\ra k^X$
sending  $e_{\sigma_i}$ to $\sum_{j\in \sigma_i}e_j$ for
$1\le i\le m$. Also, note that there is a canonical algebra homomorphism
$$\pi: kQ/\langle\rho\rangle\lra k^X$$
sending $e_i$ to $e_i$ for $i\in X$, and all other idempotents and all arrows to zero. Similarly, there is a canonical, surjective algebra homomorphism
$\pi: kQ^{\sigma}/\langle\rho^{\sigma}\rangle\ra  k^{\sigma}$.  Then we have the following commutative diagram of algebra homomorphisms:
$$\xymatrix@M=0.5mm{
kQ^{\sigma}/\langle\rho^{\sigma}\rangle  \ar[r]^{\lambda_{\sigma}}\ar[d]_{\pi} &  kQ/\langle\rho\rangle\ar[d]^{\pi}\\
*+[r]{k^{\sigma}}\ar[r]^{\lambda_{\sigma}} &*+[r]{k^X}
}$$
Moreover, $\dim_k A+\dim_k
k^{\sigma}=\dim_kA^{\sigma}+\dim_k k^X$. This implies that
the above commutative diagram is a pullback diagram.

\begin{Theo} Suppose that
$F$ is a derived equivalence between algebras $A:=kQ/\langle\rho\rangle$ and $A':=kQ'/\langle\rho'\rangle$.
Let $X$ be a subset of $Q_0$ such that the simple $A$-modules corresponding to the vertices in $X$ are sent by $F$ to simple $A'$-modules. Let $X'$ be the set of vertices in $Q_0'$ corresponding to these simple $A'$-modules.
Let $\sigma$ be a partition of $X$ and $\sigma'$ be the corresponding partition of $X'$.
Then the algebras $A^{\sigma}$ and ${A'}^{{\sigma'}}$ are derived equivalent.
\label{TheoremGluing}
\end{Theo}

{\it Proof.}
By assumption, there is a basic, radical tilting complex $\cpx{T}$
over $A$ such that $F(\cpx{T})\simeq A'$ in $\Db{A'}$. By Lemmas
\ref{LemsimpleTosimple} and \ref{LemFormofTiltingComp}, we can
rewrite $\cpx{T}$ as
   $\cpx{T}= \cpx{U}\oplus \bigoplus_{i\in X}\cpx{V_i}$
such that $\cpx{U}\in\Kb{\add(\bigoplus_{i\in Q_0\backslash X}Ae_i)}$ and $\cpx{V_i}$ is indecomposable  with $\nTP{\cpx{V_i}}{Ae_j}=\delta_{ij}$ for all $i, j\in
X$. Moreover, for each $i\in X$,
the projective $A$-module $Ae_i$ occurs as a direct summand of
$V_i^0$ with the multiplicity $1$ (see the proof of Lemma \ref{LemsimpleTosimple}).  By the definition of $\pi: A\ra k^X$, we have $k^X\otimes_AAe_i=0$ for $i\not\in X$ and  $k^X\otimes_AAe_i\simeq
k^Xe_i$ for $i\in X$. Thus there is an isomorphism in $\C{k^X}$:
   $$\cpx{h}: k^{X}\otimes_A\cpx{T}\lra k^{X}.$$
Clearly, $k^X\otimes_{k^{\sigma}}k^{\sigma}\simeq k^X$.   Let $\eta_1: \End_{\K{A}}(\cpx{T})\ra
\End_{k^X}(k^X)$ be the algebra homomorphism
determined by the composite $\pi^*\cpx{h_1}: \cpx{T}\ra
k^X\otimes_A\cpx{T}\ra k^X$, and let $\eta_2:
\End_{k^{\sigma}}(k^{\sigma})\ra
\End_{k^X}(k^X)$ be the algebra homomorphism
determined by $\lambda_{\sigma}$.  By Theorem \ref{TheoMain}, the pullback algebra of $\eta_1$ and
$\eta_2$ is derived equivalent to the pullback algebra  $A^{\sigma}$ of $\pi: A\ra k^X$ and $\lambda_{\sigma}: k^{\sigma}\ra k^{X}$. It remains to show
that ${A'}^{\sigma'}$ is isomorphic to the pullback algebra of
$\eta_1$ and $\eta_2$.

For each $x$ in $Q_0$ (respectively, $Q'_0$), we denote by
$S_x$ (respectively, $S'_x$) the  simple
$A$-module (respectively, $A'$-module) corresponding to the vertex $x$. By relabeling the vertices
if necessary, we can assume that
$$X=\{1,\cdots,m\}=X'$$ such that $F(S_i)\simeq S'_i$
for $1\le i\le m$. In this case,
 $\sigma$ and $\sigma'$ are the same partition of $\{1,\cdots,m\}$.
 For $i,j\in \{1,\cdots,m\}$, the Hom-space
  $$\Hom_{\Db{A'}}(F(\cpx{V_i}), S'_{j}))\simeq \Hom_{\Db{A'}}(F(\cpx{V_i}), F(S_j))\simeq\Hom_{\Db{A}}(\cpx{V_i}, S_j) $$
is $1$-dimensional for $i=j$, and zero for $i\neq j$. Thus it follows from the indecomposability of  $F(\cpx{V_i})$ that there exists an isomorphism $g_i: F(\cpx{V_i})\ra A'e_i$ for
$1\le i\le m$. Let $f:=\sum_{j\in Q'_0\backslash X'}e_j\in A'$.
Then there is an isomorphism $g: F(\cpx{U})\ra A'f$.  Thus we obtain an isomorphism
$${\rm diag}[g, g_1,\cdots,g_m]: F(\cpx{T})\lra A',$$
which induces an isomorphism $\tilde{g}:\End_{\Db{A'}}(F(\cpx{T}))\ra \End_{A'}(A')$. Let $s$ be the composite of the following maps
$$\End_{\K{A}}(\cpx{T})\simeq \End_{\Db{A}}(\cpx{T})\lra \End_{\Db{A'}}(F(\cpx{T}))\lraf{\tilde{g}} \End_{A'}(A')\lra A'.$$
Then, for each $i\in \{1, \cdots, m\}$, the map $s$ sends the primitive idempotent corresponding to the direct summand $\cpx{V_i}$ to $e_i$.  According to this fact, it is easy to check that the following diagram is commutative.
$$\xymatrix@M=1mm{
 \End_{\K{A}}(\cpx{T})\ar[r]^{\eta_1}\ar[d]^{\simeq}_{s} &
  \End_{k^{X}}(k^{X}) \ar[d]^{\simeq} &
  \End_{k^{\sigma}}(k^{\sigma})\ar[l]_{\eta_2}\ar[d]^{\simeq}\\
A' \ar[r]^{\pi} &
{ k^{X'}}&
  {k^{\sigma'}}\ar[l]_{\lambda_{\sigma'}}
}$$
Note that the unlabeled vertical
isomorphisms are the canonical ones. This diagram shows that ${A'}^{\sigma'}$, which is the pullback of
$\pi$ and $\lambda_{\sigma'}$, is isomorphic to the pullback algebra
of $\eta_1$ and $\eta_2$, and finishes the proof.
$\square$

\medskip
{\it Remark.} In Theorem \ref{TheoremGluing}, the indecomposable projective $A^{\sigma}$-module corresponding to a part of the partition $\sigma$ occurs only once (in degree zero) in the tilting complex that induces a derived equivalence between $A^{\sigma}$ and $A'^{\sigma'}$ (see the proof of Theorem \ref{TheoMain}). Therefore, by Lemma \ref{LemsimpleTosimple}, this derived equivalence sends the simple modules corresponding to parts of $\sigma$ to the simple modules corresponding to parts of $\sigma'$.
Thus Theorem \ref{TheoremGluing} can be applied repeatedly.

\medskip
Theorem \ref{TheoremGluing} also provides a way to construct a new derived equivalence from  two given derived equivalences.

\begin{Koro}
Let $F$ be a derived equivalence between two algebras $A:=kQ/\langle\rho\rangle$ and $A':=kQ'/\langle\rho'\rangle$, and let $G$ be a derived equivalence between $B:=k\Gamma/\langle\phi\rangle$ and $B':=k\Gamma'/\langle\phi'\rangle$. Suppose that $\bar{Q}_0$ (respectively, $\bar{\Gamma}_0$) be a subset of $Q_0$ (respectively, $\Gamma_0$) such that the simple modules corresponding to the vertices in $\bar{Q}_0$ (respectively, $\bar{\Gamma}_0$) are sent by $F$ (respectively, $G$) to simple modules corresponding to the vertices in $\bar{Q}'_0$ (respectively, $\bar{\Gamma}'_0$) and that $|\bar{Q}_0| = |\bar{Q}'_0|$ and $|\bar{\Gamma}_0| = |\bar{\Gamma}'_0|$. Let $\sigma$ be a partition of the set $\bar{Q}_0\cup\bar{\Gamma}_0$ and $\sigma'$ be the corresponding partition of $\bar{Q}'_0\cup\bar{\Gamma}'_0$. Then the algebras $(A\times B)^{\sigma}$ and $(A'\times B')^{\sigma'}$ are derived equivalent.
\label{CorGluing-1}
\end{Koro}

{\it Proof.}
Taking coproducts of algebras, we can get a derived equivalence between $A\times B$ and $A'\times B'$, which sends the simple modules corresponding to the vertices in $\bar{Q}_0\cup\bar{\Gamma}_0$ to the simple modules corresponding to the vertices in $\bar{Q}'_0\cup\bar{\Gamma}'_0$. Thus the corollary follows immediately from Theorem \ref{TheoremGluing}.
$\square$

\medskip
A special case of Corollary \ref{CorGluing-1} is $B=A^{\op}$ and $B'={A'}^{\op}$. In this case we can get derived equivalence between $(A\times A^{\op})^{\sigma}$ and $(A'\times {A'}^{\op})^{\sigma'}$ since algebras $A$ and $A'$ are derived equivalent if and only if so are their opposite algebras.

Another special case of Corollary \ref{CorGluing-1} is to attach an algebra simultaneously to derived equivalent algebras and make the resulting algebras again derived equivalent.

\begin{Koro}
Let $F$ be  a derived equivalence between the algebras $A:=kQ/\langle\rho\rangle$ and $A':=kQ'/\langle\rho'\rangle$ such that $F$ sends the simple $A$-modules corresponding to the vertices in $\bar{Q}_0$ to the simple $A'$-modules corresponding to the vertices in $\bar{Q}'_0$ and that $|\bar{Q}_0|=|\bar{Q}'_0|$. Suppose that $C:=k\Gamma/\langle\rho''\rangle$ is an arbitrary algebra. Let $\sigma$ be a partition of $\bar{Q}_0\cup \Gamma_0$ and $\sigma'$ be the corresponding partition of $\bar{Q}'_0\cup \Gamma_0$. Then the algebras $(A\times C)^{\sigma}$ and $(A'\times C)^{\sigma'}$ are derived equivalent.
\label{CorGluing-2}
\end{Koro}

\subsection{Derived equivalences from unifying arrows}

In this subsection, we shall construct new derived equivalences from given ones by unifying certain arrows in quivers.

We first fix some notation. Throughout this subsection, $\Delta$ is the quiver with the vertex set $\{x, 1, 2, \cdots, n\}$ and $n$ arrows $\alpha_j: x\ra j$, $1\le j\le n$. Here, we understand that the arrows have pairwise distinct ending vertices. We define $E:=\{1, \cdots, n\}$.  It may happen that the vertex $x$ falls into $E$. In this case $\Delta$ has the vertex set $E$. Let $\sigma$ be the partition of $E$ with only one part, and let $\alpha:=\{\alpha_1, \cdots, \alpha_n\}$ for simplicity.

\medskip
Let $A=kQ/\langle\rho\rangle$ be a finite-dimensional $k$-algebra such that $\Delta$ is a subquiver of $Q$.  By the previous discussion, there is an algebra embedding
   $$\lambda_{\sigma}: kQ^{\sigma}/\langle\rho^{\sigma}\rangle\lra kQ/\langle\rho\rangle.$$
Let $Q^{\alpha}$ be the quiver obtained from $Q^{\sigma}$ by unifying the arrows $\alpha_1, \cdots, \alpha_n$ into one arrow $\bar{\alpha}$ in $Q^{\sigma}$. Thus $Q^{\alpha}$ has the vertex set  $Q^{\sigma}$, while the set of arrows is $\{\bar{\alpha}\}\cup Q^{\sigma}_1\backslash\{\alpha_1,\cdots,\alpha_n\}$. Then there is  a canonical algebra homomorphism
$$\varphi: kQ^{\alpha}\lra kQ^{\sigma}/\langle\rho^{\sigma}\rangle$$
sending $\bar{\alpha}$ to $\sum_{i=1}^n\alpha_i$, and preserving all
other arrows and all vertices. It is easy to see that
$\Ker(\varphi)$ is contained in $\langle Q^{\alpha}_2\rangle$.  Let
$\rho^{\alpha}$ be relations on $Q^{\alpha}$ such that
$\langle\rho^{\alpha}\rangle=\Ker(\varphi)$. Then we get a natural
embedding
$$\lambda_{\alpha}: kQ^{\alpha}/\langle\rho^{\alpha}\rangle\lra kQ^{\sigma}/\langle\rho^{\sigma}\rangle.$$
We define $A^{\alpha}:=kQ^{\alpha}/\langle\rho^{\alpha}\rangle$. This is called the \emph{unifying algebra} of $A$ by $\alpha$. The
image of the composite $\lambda_{\alpha}\lambda_{\sigma}$
is the subalgebra of $A$ generated by all the arrows
$\beta\not\in\alpha$, $\sum_{i=1}^n\alpha_i$ and idempotents $e_{E}$, $e_i, i\in Q_0\backslash E$. The above procedure can be illustrated visually by the following (local)
pictures: 
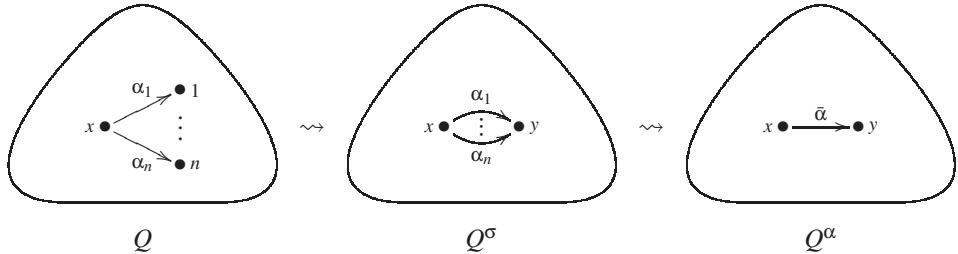
\vspace{-0.3cm}
$$\xy
    (0,0)*+{\bullet}="a", (10, 5)*+{\bullet}="1", (10, -5)*+{\bullet}="n",
    {\ar "a"; "1"},
    {\ar "a"; "n"},
    (10, 1)*+{\vdots},
    (5,5)*+{\scriptstyle{\alpha_1}},
    (5,-5)*+{\scriptstyle{\alpha_n}},
    (12,5)*+{\scriptstyle{1}},
    (12,-5)*+{\scriptstyle{n}},
    (-2,0)*+{\scriptstyle{x}},
    (-5,-10)*+0{}="l", (15,-10)*+0{}="r",
    "l";"r" **\crv{(-20,-10)&(0,18)&(10,18)&(30,-10)},
    **\dir{-},
    (5, -15)*+{Q}
  \endxy \rightsquigarrow \xy
    (0,0)*+{\bullet}="a", (10, 0)*+{\bullet}="y",
    (5,1)*+{\cdot}, (5,0)*+{\cdot},
    (5,-1)*+{\cdot},
    {\ar@/^/^{\alpha_1} "a"; "y"},
    {\ar@/_/_{\alpha_n} "a"; "y"},
    (12,0)*+{\scriptstyle{y}},
    (-2,0)*+{\scriptstyle{x}},
    (-5,-10)*+0{}="l", (15,-10)*+0{}="r",
    "l";"r" **\crv{(-20,-10)&(0,18)&(10,18)&(30,-10)},
    **\dir{-},
    (5, -15)*+{Q^{\sigma}}
  \endxy\rightsquigarrow \xy
    (0,0)*+{\bullet}="a", (10, 0)*+{\bullet}="y",
    {\ar^{\bar{\alpha}} "a"; "y"},
    (12,0)*+{\scriptstyle{y}},
    (-2,0)*+{\scriptstyle{x}},
    (-5,-10)*+0{}="l", (15,-10)*+0{}="r",
    "l";"r" **\crv{(-20,-10)&(0,18)&(10,18)&(30,-10)},
    **\dir{-},
    (5, -15)*+{Q^{\alpha}}
  \endxy$$
Next, we shall interpret the
algebra $A^{\alpha}$ as a pullback algebra.  Actually, $A^{\alpha}$ fits into the following pullback diagram of algebra homomorphisms:
$$\xymatrix@M=1mm{
A^{\alpha} \ar[r]^{\lambda_{\alpha}}\ar[d]_(.4){\pi} & A^{\sigma}\ar[d]^(.45){\pi}\\
\quad k^{\Delta^{\sigma}_0}\quad \ar[r]^(.35){\lambda} &
\quad (k\Delta)^{\sigma}/\langle\sum_{i=1}^n\alpha_i\rangle
}$$
The vertical homomorphisms in the above diagram are
obviously defined.

\begin{Lem}
The algebra $(k{\Delta})^{\sigma}/\langle\sum_{i=1}^n\alpha_i\rangle$ is radical-square zero.
   \label{LemmaRadicalsquarezero}
\end{Lem}

{\it Proof.}
If $x\not\in\{1,\cdots, n\}$, then $x\neq y$ and
$(k{\Delta})^{\sigma}/\langle\sum_{i=1}^n\alpha_i\rangle$
is radical-square zero. Without loss of generality, we now assume
that $\alpha_1$ is a loop in the quiver $\Delta$. Then none of
$\alpha_2,\cdots,\alpha_n$ is a loop by the assumption that the
vertices $1,\cdots,n$ are pairwise distinct. Thus,
$\alpha_i\alpha_j=0$ for all $i\neq 1$ and all $j\in\{1,\cdots,n\}$.
Further, for each $j\in \{1,\cdots, n\}$, the path $\alpha_1\alpha_j=
(\sum_{i=1}^n\alpha_i)\alpha_j$ is in
$\langle\sum_{i=1}^n\alpha_i\rangle$. Altogether, we have shown that all paths in
$(k\Delta)^{\sigma}$  of length $2$ belong to
$\langle\sum_{i=1}^n\alpha_i\rangle$, and the lemma is then proved.
$\square$

\medskip
Let $kQ/\langle\rho\rangle$ be a finite-dimensional algebra defined by
a quiver $Q$ with relations $\rho$. Let $i$ and $j$ be vertices in $Q_0$,
and let $Q_{ij}$ be the $k$-vector space with all arrows from
$i$ to $j$ as a basis. Then every vector space automorphism  $\chi:
Q_{ij}\ra Q_{ij}$ extends to an algebra automorphism $\phi_{\chi}:
kQ\ra kQ$ which sends $\alpha\in Q_{ij}$ to
$(\alpha)\chi$ and preserves all other arrows and all vertices. If
$(\langle\rho\rangle)\phi_{\chi}=\langle\rho\rangle$ for all such
automorphisms $\chi$ on $Q_{ij}$, then $\rho$ is said to be \emph{
$(i,j)$-invariant}. Let $\Gamma=(\Gamma_0,\Gamma_1)$ be a sub-quiver of $Q$. We say that $\rho$ is  \emph{$\Gamma$-invariant} if $\rho$ is $(i,j)$-invariant for all $i, j\in \Gamma_0$. For example, $\rho$ is $\Gamma$-invariant if $\rho$ consists only of monomial relations and there is at most $1$ arrow from $i$ to $j$ in $Q$ for any two vertices $i,j$ in $\Gamma_0$. Note that $\rho$ is $\Gamma$-invariant if and only if $\rho^{\opp}$ in $kQ^{\opp}/\langle\rho^{\opp}\rangle$ is $\Gamma^{\opp}$-invariant.

\begin{Theo}
Let $A:=kQ/\langle\rho\rangle$ and
$A':=kQ'/\langle\rho'\rangle$ be algebras, and suppose that the given quiver $\Delta$ is a sub-quiver of both $Q$ and $Q'$.
Assume that $\rho$ or $\rho'$ is $\Delta$-invariant.
If $F:\Db{A}\ra \Db{A'}$ is a derived equivalence such that $F(S_i)\simeq S'_i$ for all $i\in \Delta_0$, then $A^{\alpha}$ and ${A'}^{\alpha}$ are
derived equivalent.
 \label{TheoUnifyingArrows}
\end{Theo}

{\it Proof.}
Without loss of generality, we assume that $\rho'$ is $\Delta$-invariant.
Further, we assume that the common starting
vertex $x$ of $\alpha_1,\cdots,\alpha_n$ is not in $E$. The
case that $x\in E$ can be proved similarly.  Let $\tilde{\Delta}$ be the full sub-quiver of $Q$ defined by  $\Delta_0$. Then $\Delta$ is a sub-quiver of $\tilde{\Delta}$ with the same vertices and (possibly) less arrows. Let $B:=k\tilde{\Delta}/\langle\tilde{\Delta}_2\rangle$, and let $\Lambda:=(k\Delta)^{\sigma}/\langle\sum_{i=1}^n\alpha_i\rangle$.  Then, by Lemma \ref{LemmaRadicalsquarezero}, there is a canonical surjective homomorphism $\pi: B^{\sigma}\ra \Lambda$ of algebras.

Let $\cpx{T}$ be a basic, radical
tilting complex associated to the derived equivalence $F$. Set $U:=\bigoplus_{i\in
Q_0\backslash\Delta_0}Ae_i.$
%$\displaystyle{U:=\bigoplus_{i\in
%Q_0\backslash\Delta_0}Ae_i}.$
Since $F(S_i)\simeq S'_i$ for
all $i\in\Delta_0$, we can assume
$$\cpx{T}=\cpx{U}\oplus\cpx{V_x}\oplus\cpx{V_1}\oplus\cdots\oplus\cpx{V_n}$$
by Lemmas \ref{LemsimpleTosimple} and
\ref{LemFormofTiltingComp}, where $\cpx{V_i}$ is a complex in $\Kb{\pmodcat{A}}$ such that, for
each $i\in\Delta_0$,  $V_i^0=Ae_i\oplus U_i$ for some
$U_i\in\add(U)$ and $V_i^j\in\add(U)$ for all $j\neq 0$.
Note that there is a commutative diagram
$$\xymatrix@M=0.5mm{
A^{\sigma} \ar[r]^{\pi}\ar[d]^{\lambda_{\sigma}} & B^{\sigma}\ar[d]^{\lambda_{\sigma}}\ar[r]^{\pi} & k^{\sigma}\ar[d]^{\lambda_{\sigma}} \\
A\ar[r]^{\pi} & B \ar[r]^{\pi} & k^E\\
}$$
where the horizontal maps are the canonical maps. The right-hand square and the entire square are pullback diagrams of algebras. This implies that the left-hand square is also a pullback diagram.  It is easy to see  $B\otimes_AU=0$ and that there is an isomorphism of stalk complexes in $\C{B}$:
$$\cpx{h}: B\otimes_A\cpx{T}=B\otimes_A(Ae_x\oplus \bigoplus_{i=1}^nAe_i)\lra B\simeq B\otimes_{B^{\sigma}}B^{\sigma}.$$
By the proof of Theorem \ref{TheoMain},  the complex $\cpx{T_{\sigma}}:=M(\cpx{T}, B^{\sigma}, \cpx{h})$
is a tilting complex over $A^{\sigma}$ with $\End_{\K{A^{\sigma}}}(\cpx{T_{\sigma}})\simeq {A'}^{\sigma}$.  Moreover, there is a pullback diagram
$$\vcenter{\xymatrix@M=1mm{
\End_{\K{A^{\sigma}}}(\cpx{T_{\sigma}})\ar[r]^{\epsilon_1}\ar[d]^{\epsilon_2\mu} &\End_{\K{A}}(\cpx{T})\ar[d]^{\eta\mu}\\
B^{\sigma}\ar[r]^{\lambda_{\sigma}} &B,
}}$$
where $\eta$ is determined by $\cpx{T}\ra B$, $\epsilon_1$ and $\epsilon_2$ are determined by the projections from $\cpx{T_{\sigma}}$ to $\cpx{T}$ and $B^{\sigma}$, respectively, and $\mu$ is the canonical isomorphism from $\End(_{\Lambda}\Lambda)$ to $\Lambda$ for an algebra $\Lambda$.

\medskip
By assumption, $F(S_i)\simeq S_i'$ for all $i\in\Delta_0$. It follows that $\Ext_A^1(S_i, S_j)\simeq\Ext_{A'}^1(S_i', S_j')$ for all $i, j\in\Delta_0$. This indicates that the number of arrows from $i$ to $j$ are equal in both $Q$ and $Q'$. Hence we can assume that $\tilde{\Delta}$ is also a full sub-quiver of $Q'$ with vertices $\Delta_0$.  As a consequence, there is a canonical, surjective homomorphism $\pi: A'\ra B$ of algebras.

Let $\theta: A'\ra \End_{\K{A}}(\cpx{T})$ be an isomorphism of algebras.  Note that $\End_B(B)\simeq B$ is radical-square zero by definition. Thus it is easy to know that $\theta \eta\mu: A'\ra B$ sends the kernel of $\pi: A'\ra B$ to zero, and that there is an algebra homomorphism $\chi: B\ra B$,  which fixes all idempotents $e_i, i\in\Delta_0$, such that $\theta \eta\mu=\pi\chi$.  Since $\chi$ fixes the idempotents $e_i$ with $i\in\Delta_0$, it induces an automorphism of the vector space $e_iBe_j$ which is isomorphic to the vector space $Q'_{ij}$ for all $i, j\in \Delta_0$. Since $\rho'$ is $\Delta$-invariant, there is an automorphism $\phi_{\chi}: A'\ra A'$ extending $\chi$, that is, $\phi_{\chi}\pi=\pi\chi$. Thus $\theta^{-1}\phi_{\chi} \pi=\eta\mu$, that is,  there is a commutative diagram
$$\xymatrix@M=1mm{
\End_{\K{A}}(\cpx{T})\ar[r]^(.65){\eta\mu}\ar[d]^{\theta^{-1}\phi_{\chi}}_{\simeq} & B\ar@{=}[d] &\ar[l]_{\lambda_{\sigma}} B^{\sigma}\ar@{=}[d]\\
A'\ar[r]^{\pi} & B &B^{\sigma}\ar[l]_{\lambda_{\sigma}}\\
}$$
It then follows that there is an isomorphism $\psi$ from the pullback algebra $\End_{\K{A^{\sigma}}}(\cpx{T_{\sigma}})$ of $\eta\mu$ and $\lambda_{\sigma}$ to the pullback algebra ${A'}^{\sigma}$ of $\pi$ and $\lambda_{\sigma}$ such that  the following diagram
$$\xymatrix@M=1mm{
\End_{\K{A^{\sigma}}}(\cpx{T_{\sigma}})\ar[r]^{\epsilon_2}\ar[d]^{\psi} & \End_{B^{\sigma}}(B^{\sigma})\ar[d]^{\mu}\\
{A'}^{\sigma}\ar[r]^{\pi} & B^{\sigma}
}$$
is commutative. This diagram can be extended to the following commutative diagram
$$\xymatrix@M=1.5mm{
\End_{\K{A^{\sigma}}}(\cpx{T_{\sigma}})\ar[r]^{\epsilon_2}\ar[d]^{\psi}_{\simeq} & \End_{B^{\sigma}}(B^{\sigma})\ar[d]^{\mu}_{\simeq} \ar[r]^{p} & \End_{\Lambda}(\Lambda)\ar[d]^{\mu}_{\simeq} & \End_{k^{\Delta_0^{\sigma}}}(k^{\Delta_0^{\sigma}})\ar[l]_{i}\ar[d]^{\mu}_{\simeq}\\
{A'}^{\sigma}\ar[r]^{\pi} & B^{\sigma}\ar[r]^{\pi} &\Lambda & k^{\Delta_0^{\sigma}}\ar[l]_{\lambda}
}$$
where $p$ and $i$ are determined by $\pi$ and $\lambda$, respectively. It then follows that the pullback algebra ${A'}^{\alpha}$ of $\pi: {A'}^{\sigma}\ra \Lambda$ and $\lambda$ is isomorphic to the pullback algebra of $\epsilon_2p$ and $i$.  Note that  $$\Lambda\otimes_{k^{\Delta_0^{\sigma}}}k^{\Delta_0^{\sigma}}\simeq \Lambda \simeq \Lambda\otimes_{B^{\sigma}}B^{\sigma}\simeq \Lambda\otimes_{B^{\sigma}}B^{\sigma}\otimes_{A^{\sigma}}\cpx{T_{\sigma}}$$
in $\C{\Lambda}$. By the proof of Theorem \ref{TheoMain}, the pullback algebra of $\epsilon_2p$ and $i$ is derived equivalent to the pullback algebra $A^{\alpha}$ of $A^{\sigma}\lraf{\pi} \Lambda\llaf{\lambda} k^{\Delta_0^{\sigma}}$. Consequently, ${A'}^{\alpha}$ is derived equivalent to $A^{\alpha}$.
$\square$

\begin{Rem} \rm{ (1) Note that two algebras $A$ and $B$ are derived equivalent if and only if their opposite algebras $A^{\opp}$ and $B^{\opp}$ are derived equivalent. So we can replace $\Delta$ by $\Delta^{\opp}$ and consider unifying arrows of $\Delta^{\opp}$. This means that Theorem \ref{TheoUnifyingArrows} also holds true for the subquiver $\Delta^{\opp}$.}

(2)
The derived equivalence constructed in theorem \ref{TheoUnifyingArrows} sends the simple $A^{\alpha}$-modules corresponding to $x$ and $y$ again to simple $A'^{\alpha}$-modules.
\end{Rem}

\subsection{Derived equivalences from identifying socle elements}\label{subsection-glue-e-cycles}

In this subsection, we introduce the third operation by identifying socle elements of algebras to get new derived equivalences.

Let $A$ be a basic Artin algebra with the Jacobson radical $r_A$, and let $1_A=e_1+\cdots +e_n$ be a decomposition of $1_A$ into pairwise orthogonal primitive idempotents. Fix $i, j\in \{1, \cdots, n\}$. A {\em longest $(e_i, e_j)$-element}  in $A$ is a nonzero element $a\in e_ir_Ae_j$ such that $r_Aa=0=ar_A$, that is, $a\in \soc(r_Ae_j)\cap \soc(e_ir_A)$. In this case, the ideal $\langle a\rangle$ of $A$ generated by $a$ is $1$-dimensional and contained in $\soc({}_AAe_j)\cap\soc(e_iA{}_A)$.  A longest  $(e_i, e_i)$-element is called a {\em complete $e_i$-cycle}.

\medskip
For the rest of this subsection, we fix two algebras $A:=kQ/\langle\rho\rangle$ and $B:=k\Gamma/\langle\omega\rangle$ given by quivers with relations. Suppose that  $a$ is a longest $(e_i, e_j)$-element in $A$, and that  $b$ is a longest $(e_{s}, e_{t})$-element in $B$,  where $i, j\in Q_0$ and $s, t\in \Gamma_0$. We glue $i$ and $s$ into a new vertex $u$, and glue $j$ and $t$ into a new vertex $v$. Let $\sigma$ be the corresponding partition of the set $\{i, j, s, t\}$. In case that $i=j$ or $s=t$, we actually glue all the vertices into one vertex, that is, $u=v$. Let $(A\times B)^{\sigma}$ be the $\sigma$-gluing algebra defined in Subsection \ref{subsection-gluing-vertices}. In case that $i=j$ and $s=t$, we simply write $A\glueidem{e_i}{e_s}B$ for $(A\times B)^{\sigma}$.  Now, it is easy to see that $a-b$ is a longest $(e_u, e_v)$-element in $(A\times B)^{\sigma}$ and the ideal $\langle a-b\rangle$ of $(A\times B)^{\sigma}$ generated by $a-b$ is $1$-dimensional. So, we can define a new algebra
 $$A\gluesocab B:=(A\times B)^{\sigma}/\langle a-b \rangle.$$
It is called the algebra of \emph{identifying socle elements} in $A$ and $B$.

Suppose that $A':=kQ'/\langle\rho'\rangle$ is another algebra and there is a derived equivalence $F:\Db{A}\ra \Db{A'}$ such that $F(S_i)\simeq S_{i'}$ and $F(S_j)\simeq S_{j'}$ for some $i', j'\in Q'_0$.  Let $\cpx{T}$ be a basic, radical tilting complex associated to $F$. We may identify $A'$ with $\End_{\Kb{A}}(\cpx{T})$ via the isomorphism $\End_{\Kb{A}}(\cpx{T})\ra A'$ induced by $F$. Further, by the proof of Lemma \ref{LemsimpleTosimple}, both $Ae_i$ and $Ae_j$ only occur in degree zero with  the multiplicity $1$ in $\cpx{T}$. For $x\in \{i, j\}$, let $\cpx{T_x}$ be  the indecomposable direct summand of $\cpx{T}$ such that $Ae_x$ is a direct summand of $T_x^0$, namely $T_x^0=Ae_x\oplus P_x$, and let $e_{x'}$ be the primitive idempotent element in $A'$ corresponding to the summand $\cpx{T_x}$.   Let  $m_a: \cpx{T_i}\ra \cpx{T_j}$ be the following (well-defined) particular morphism
$$\xymatrix@M=1mm{
\cdots\ar[r] & T_i^{-1}\ar[d]^{0}\ar[r] & Ae_i\oplus P_i\ar[r]\ar[d]^{\left[\begin{smallmatrix}\cdot a & 0\\0&0\end{smallmatrix}\right]} & T_i^1\ar[r]\ar[d]^0 &\cdots\\
\cdots\ar[r] & T_j^{-1}\ar[r] & Ae_j\oplus P_j\ar[r] & T_j^1\ar[r] &\cdots\, \
}$$
and let ${a'}$ be the composite $\cpx{T}\ra\cpx{T_i}\lraf{m_a}\cpx{T_j}\ra \cpx{T}$, where the first and last morphisms are the canonical projection and injection, respectively. This element $a'$ has the following property.

\begin{Lem} The element $a'$ just defined is a longest $({e}_{i'}, {e}_{j'})$-element in $A'$.\label{completecycle}
\end{Lem}

{\it Proof.} Since $a\in e_ir_Ae_j$ is nilpotent, the element $a'$ is nilpotent and lies in ${e}_{i'}r_{A'}{e}_{j'}$. It remains to show $r_{A'}a'=0=a'r_{A'}$.

Let $\cpx{g}: \cpx{T}\ra \cpx{T}$ be in $r_{A'}$. Then $\cpx{g}$ is nilpotent, that is, $(\cpx{g})^m$ is null-homotopic for some integer $m\ge 1$. Particularly, $(g^0)^m=h^0d^{-1}+d^0h^1$ for some  homomorphisms $h^0: T^0\ra T^{-1}$ and $h^1: T^1\ra T^{0}$ of $A$-modules.
Since the differential maps of $\cpx{T}$ are radical by assumption, the map $(g^0)^m$ is radical, and so is $g^0$. It follows that the composite $g^0\pi^0: T^0\ra T_i^0$ is also a radical map, where $\cpx{\pi}$ is a canonical projection $\cpx{T}\ra \cpx{T_i}$.  Now, the fact $r_Aa=0$ indicates that the composite $Ae_l\raf{r} Ae_i\raf{\cdot a} Ae_j$ is zero for all $l\in Q_0$ and all radical maps $r$.  Hence the chain map $\cpx{g}\cpx{\pi}m_a$ is zero in all degrees, and consequently $\cpx{g}a'=0$.  This shows $r_{A'}a'=0$. Similarly, using $ar_A=0$, we can prove $a'r_{A'}=0$.  $\square$

The following theorem shows that we can extend the derived equivalence between $A$ and $A'$ by  identifying socle elements.

\begin{Theo}
The algebras
$A\gluesocab B$ and ${A'}\gluesoc{{a'}}{b}B$ are derived equivalent.  \label{theorem-gluesoc}
\end{Theo}

{\it Proof.}
For simplicity, we write $\Lambda$ for $(A\times B)^{\sigma}$. As explained in Subsection \ref{subsection-gluing-vertices}, $\Lambda$ is the pullback algebra of the canonical surjective homomorphisms $B\ra k^{\sigma}$ and $A\ra k^{\sigma}$. Let $\sigma'=\{i',s\}\cup\{j',t\}$ be the corresponding partition of $\{i', j', s, t\}$.  By the proof of Theorem \ref{TheoMain}, the complex $\cpx{\tilde{T}}:=M(\cpx{T}, B, 1)$ is a tilting complex over $\Lambda$ with the endomorphism algebra isomorphic to $(A'\times B)^{\sigma'}$.  By definition, $\cpx{\tilde{T}_i}:=M(\cpx{T_i}, Be_s, 1)$ and $\cpx{\tilde{T_j}}:=M(\cpx{T_j}, B_t, 1)$ are indecomposable direct summands of $\cpx{\tilde{T}}$. Note that all other indecomposable direct summands of $\cpx{\tilde{T}}$ are of the form $M(\cpx{P}, 0, 0)$ or $M(0, Q, 0)$, where $\cpx{P}$ is an indecomposable direct summand of $\cpx{T}$ and $Q$ is an indecomposable projective $B$-module. Moreover, the indecomposable projective $\Lambda$-modules $\Lambda e_u$ and $\Lambda e_v$, which are isomorphic to $M(Ae_i, Be_s, 1)$ and $M(Ae_j, Be_t, 1)$, respectively, only occur in degree zero with the multiplicity $1$ in $\cpx{\tilde{T}}$. Thus $\cpx{\tilde{T}}$ is a basic, radical complex over $\Lambda$.

Set $I:=\langle a-b\rangle$. Then $Ie_v=I=e_uI$ and $IX=0$ for all indecomposable projective $\Lambda$-modules $X$ not isomorphic to $\Lambda e_v$.
It follows that $I\cpx{\tilde{T}}=I \tilde{T}^0\simeq {}_{\Lambda}I$. Note that ${}_{\Lambda}I$ is a simple $\Lambda$-module with $e_uI\neq 0$. Hence $\Hom_{\Kb{\Lambda}}(\cpx{\tilde{T}}, I\cpx{\tilde{T}}[l])\simeq \Hom_{\Kb{A}}(\cpx{\tilde{T}}, I[l])=0$ for all $l\neq 0$. Now, the short exact sequence
$0\ra I\cpx{\tilde{T}}\ra\cpx{\tilde{T}}\ra\cpx{\tilde{T}}/I\cpx{\tilde{T}}\ra 0$
in the category of complexes over $\Lambda$ gives raise to a triangle
$$I\cpx{\tilde{T}}\lra\cpx{\tilde{T}}\lra\cpx{\tilde{T}}/I\cpx{\tilde{T}}\lra I\cpx{\tilde{T}}[1]$$
in $\Db{\Lambda}$. Applying $\Hom_{\Db{\Lambda}}(\cpx{\tilde{T}}, -)$ to this triangle, we get an exact sequence
$$0\lra \Hom_{\Db{\Lambda}}(\cpx{\tilde{T}}, \cpx{\tilde{T}}/I\cpx{\tilde{T}}[-1])\lra \Hom_{\Db{\Lambda}}(\cpx{\tilde{T}}, I\cpx{\tilde{T}})\lra \Hom_{\Db{\Lambda}}(\cpx{\tilde{T}}, \cpx{\tilde{T}}), $$
which is isomorphic to
$$ (\sharp)\qquad 0\lra \Hom_{\Kb{\Lambda}}(\cpx{\tilde{T}}, \cpx{\tilde{T}}/I\cpx{\tilde{T}}[-1])\lra \Hom_{\Kb{\Lambda}}(\cpx{\tilde{T}}, I\cpx{\tilde{T}})\lraf{\theta} \Hom_{\Kb{\Lambda}}(\cpx{\tilde{T}}, \cpx{\tilde{T}}).  $$
Note that the map $\cdot (a-b): \Lambda e_u\ra {}_{\Lambda}I$ induces a morphism $\cpx{g}$ in $\End_{\Kb{\Lambda}}(\cpx{\tilde{T}_i})$:
$$\xymatrix@M=0.2mm{
\cdots \ar[r] & T_i^{-1}\ar[r]^(.35){d}\ar[d] & \Lambda e_u\oplus P_i\ar[r]\ar[d] ^{\left[\begin{smallmatrix} \cdot (a-b)\\ 0\end{smallmatrix}\right]}& T_i^{1} \ar[d]\ar[r] &0\\
&\; 0 \; \ar[r]\ar[d]& \; I \; \ar[r]\ar[d] & \; 0 \; \ar[d]\\
0\ar[r] & T_j^{-1}\ar[r]^(.35){d} & \Lambda e_v\oplus P_j\ar[r] & T_j^1 \ar[r] &0.\\
}$$
The image of $g^0$ is $I\Lambda e_v=I$. It follows that $\cpx{g}$ cannot be null-homotopic, since the image of  any morphism from $T_j^{-1}$ or $T_i^1$ to $\Lambda e_v$ has image contained in $Ae_j$ which intersects $I$ trivially.  Hence $\cpx{g}\neq 0$, and therefore
$$\cpx{\tilde{g}}:=\begin{bmatrix}\cpx{g} &0\\0&0\end{bmatrix}$$
is a nonzero endomorphism of $\cpx{\tilde{T}}$ and lies in $\Img(\theta)$ (see the sequence $(\sharp)$). Note that $$\Hom_{\Kb{\Lambda}}(\cpx{\tilde{T}}, I\cpx{\tilde{T}})\simeq \Hom_{\Kb{\Lambda}}(\cpx{\tilde{T}}, {}_{\Lambda}I)\simeq \Hom_{\Lambda}(\Lambda e_u, {}_{\Lambda}I)\simeq e_uI=I$$
and $I$ is $1$-dimensional. Hence $\theta$ is an injective map and $\Img(\theta)$ is a $1$-dimensional $k$-space with $\tilde{g}$ as a basis. It follows from ($\sharp$) that
$\Hom_{\Kb{\Lambda}}(\cpx{\tilde{T}}, \cpx{\tilde{T}}/I\cpx{\tilde{T}}[-1])=0$. Thus
$\Hom_{\Kb{\Lambda}}(\cpx{\tilde{T}}/I\cpx{\tilde{T}}, \cpx{\tilde{T}}/I\cpx{\tilde{T}}[-1])\simeq \Hom_{\Kb{\Lambda}}(\cpx{\tilde{T}}, \cpx{\tilde{T}}/I\cpx{\tilde{T}}[-1])=0.$
Now, by \cite[Theorem 4.2]{HuXi2013}, the algebras $\Lambda/I$ and $\End_{\Kb{\Lambda}}(\cpx{\tilde{T}})/\Img(\theta)$ are derived equivalent. It is easy to check that the isomorphism $\End_{\Kb{\Lambda}}(\cpx{\tilde{T}})\simeq (A'\times B)^{\sigma'}$, which is induced by the projections $\Lambda\ra A$ and $\Lambda\ra B$, sends the element $\cpx{\tilde{g}}$ in $\End_{\Kb{\Lambda}}(\cpx{\tilde{T}})$ to  $a'-b$ in $(A'\times B)^{\sigma'}$. As a result, $\End_{\Kb{\Lambda}}(\cpx{\tilde{T}})/\Img(\theta)$   is isomorphic to $A'\gluesoc{a'}{b}B$. Note that the algebra $\Lambda/I$ is just $A\gluesocab B$. Hence $A\gluesocab B$ is derived equivalent to $A'\gluesoc{a'}{b}B$. This finishes the proof.
$\square$

Remark that the derived equivalence in Theorem \ref{theorem-gluesoc} sends the simple modules over $A\gluesocab B$ corresponding to $u$ and $v$ also to simple modules over $A'\gluesoc{a'}{b}B$ corresponding to $u'$ and $v'$.

A special case of Theorem \ref{theorem-gluesoc} is that we take complete cycles with the same starting and ending vertices.

\begin{Koro}
Suppose that $e$ and $f$ are primitive idempotent elements in $A$ and $B$, respectively, and that $a\in A$ is a complete $e$-cycle and $b\in B$ is a complete $f$-cycle. Let $\cpx{T}$ be a basic, radical tilting complex over $A$ with $\nTP{\cpx{T}}{Ae}=1$, and let $A' =\End_{\Kb{A}}(\cpx{T})$. Then
$A\gluesocab B$ and $A'\gluesoc{a'}{b}B$ are derived equivalent.  \label{cor-identifysoc}
\end{Koro}

\section{Derived equivalences and Frobenius type \label{sect5} }

As an application of our constructions in Section \ref{sect4}, we consider, in this section, whether Frobenius type of algebras is invariant under derived equivalences. Solutions to this question are presented in Proposition \ref{tiltinginvariant}, Corollary \ref{cor3} and Example \ref{Bsp5.6}.

Throughout this section, all algebras are finite-dimensional over a field.

Frobenius parts of algebras have played an important role in several aspects of the representation theory of algebras. For instance, concerning the Auslander-Reiten conjecture (or Alperin-Auslander conjecture referred in \cite{Rouquier2006}) which states that stable equivalent algebras should have the same number of non-isomorphic, non-projective simple modules, Mart\'inez-Villa reduced the validity of this conjecture for algebras without nodes to that for Frobenius parts ({see \cite{Martinez-Villa1990a}). In \cite{HuXi2014-preprint}, the problem of lifting stable equivalences of Morita type to derived equivalences for arbitrary algebras is reduced to the one for their Frobenius parts. Moreover, there are close connections between dominant dimensions and Frobenius parts of algebras (see \cite{Chen2016}).

Let $A$ be an algebra. We may suppose $\stp{A}=\add(Ae)$ for $e$ an idempotent element in $A$ and that $Ae$ is a basic $A$-module. Following \cite{HuXi2014-preprint}, the algebra $eAe$ is called the \emph{Frobenius part} of $A$. It is a self-injective algebra introduced first in \cite{Martinez-Villa1990a} (see also \cite[Lemma 2.5]{HuXi2014-preprint}) and uniquely determined by $A$ up to Morita equivalence.
We say that
$A$ is \emph{Frobenius-finite (-tame or -wild)} if its Frobenius part is representation-finite (-tame or -wild). By Frobenius type we mean the representation type of the Frobenius part. If the Frobenius part of $A$ is zero, we say that $A$ is \emph{Frobenius-free}.

Frobenius-finite algebras include representation-finite algebras, Auslander-algebras, cluster-tilted algebras, and can be produced from triangular matrix algebras, Auslander-Yoneda algebras and Frobenius extensions (see \cite[Section 5]{HuXi2014-preprint} for details). For Frobenius-finite algebras over an algebraically closed field, every stable equivalence of Morita type lifts to a derived equivalence (see \cite[Theorem 1.1]{HuXi2014-preprint}). Thus this large class of algebras shares many stable and derived invariants (see \cite{Rickard1989a, Rickard1989, Krause1998, Martinez-Villa1990a, HuXi2013})

\medskip
Now, we consider behaviors of Frobenius type under stable and derived equivalences.

From \cite[Lemma 3.3]{HuXi2014-preprint} it follows that, for indecomposable algebras with separable semisimple quotients by their Jacobson radicals, Frobenius type is preserved by stable equivalences of Morita type.
Actually, a more general statement is true, namely stable equivalences preserve Frobenius type. This follows from a simple observation (see Proposition \ref{tiltinginvariant}(1) below). Recall that a simple $A$-module $S$ is called a \emph{node} in \cite{Martinez-Villa1980a} if it is neither projective, nor injective, and the almost split sequence $0\ra S \ra P\ra \mbox{Tr}D(S) \ra 0$ has a projective middle term $P$. For the definition of almost split sequences, we refer the reader to \cite{Auslander1995}.

\begin{Prop} $(1)$ Let $A$ and $B$ be algebras over an algebraically closed field and without nodes. If $A$ and $B$ are stably equivalent (that is, the stable categories $\stmodcat{A}$ and $\stmodcat{B}$ are equivalent), then they have the same Frobenius type.

$(2)$ Let $A$ be an algebra over an arbitrary field and $_AT$ be a tilting $A$-module with $B:=\End_A(T)$. Then the Frobenius parts of $A$ and $B$ are isomorphic.
\label{tiltinginvariant}
\end{Prop}

{\it Proof.}
(1) Under the assumptions of the proposition, we know from \cite{Martinez-Villa1990a} that a stable equivalence between $A$ and $B$ induces a stable equivalence between their Frobenius parts. Since stable equivalences preserve representation type by \cite{Krause1998}, we see that the Frobenius parts of $A$ and $B$ have the same representation type, and therefore $A$ and $B$ have the same Frobenius type.

(2) This follows from \cite[Lemma 4.3]{Chen2016}.
$\square$

As is known, derived equivalences between self-injective algebras over a field preserve representation type. Also, by Proposition \ref{tiltinginvariant}(2), derived equivalences induced by tilting modules over arbitrary algebras preserve Frobenius type. Furthermore, almost $\nu$-stable derived equivalences also preserve Frobenius type (see \cite[Proposition 3.3]{HuXi2014-preprint}). So, based on these phenomena, one may naturally ask the following question:

\medskip
{\bf Question.} Does a derived equivalence always preserve Frobenius type of algebras?

\medskip
In the following, we shall answer the question negatively.

Let $A$ and $B$ be basic algebras, and let $e$ and $f$ be primitive idempotents in $A$ and $B$, respectively. Suppose that $a\in A$ is a complete $e$-cycle and that $b\in B$ is a complete $f$-cycle. Set $\Lambda:=A\glueidemef B$, and $\Gamma:=A\gluesocab B$.  Recall that $\Gamma$ is the quotient algebra of $\Lambda$ modulo the one-dimensional ideal $I:=\langle a-b\rangle$.
For $x\in \Lambda$, we write $\bar{x}= x+I$ in $\Gamma$. As before, let $1_A=e+e_2+\cdots+e_n$ and $1_B=f+f_2+\cdots+f_m$ be decompositions of identities into pairwise orthogonal primitive idempotents.
Then $1_{\Gamma}= \bar{e+f}+\bar{e_2}+\cdots +\bar{e_n}+\bar{f_2}+\cdots +\bar{f_m}$ is a decompositions of $1_{\Gamma}$ into pairwise orthogonal primitive idempotents.

\medskip
In the following, we describe the Frobenius part of the algebra $\Gamma$.

\begin{Lem}
Let $A$ be an algebra, and let $e_1, e_2$ be primitive idempotents in $A$. Then the following are equivalent.

$(1)$ $\nu_A(Ae_1)\simeq Ae_2$.

$(2)$  $e_1\soc (Ae_2)\neq 0$ and, for each   $0\neq u\in e_1\soc (Ae_2)$, the following two conditions are satisfied:

\quad\,\, ${\rm (i)}$ For each $0\neq x\in Ae_2$, there is an element $y\in e_1A$ such that $yx=u$.

\quad\,\, ${\rm (ii)}$  For each $0\neq y\in e_1A$, there is an element $x\in Ae_2$ such that $yx=u$.

$(3)$ There is a nonzero element $u\in e_1Ae_2$, satisfying the conditions {\rm (i)} and ${\rm (ii)}$ in $(2)$.
\label{lemma-Ae-Nakayama}
\end{Lem}

{\it Proof.}
$(1)\Rightarrow (2)$   Suppose $\nu_A(Ae_1)\simeq Ae_2$. Then $\soc(Ae_2)$ is isomorphic to the top of $Ae_1$. Hence $e_1\soc(Ae_2)\neq 0$. Let $u$ be a nonzero element in $e_1\soc(Ae_2)$. We claim that $u\in \soc(e_1A)$.  For $r\in r_A$, let $\phi_r: Ae_2\ra A, z\mapsto zr$.
Then $\phi_r$ is a homomorphism of left $A$-modules. Since $\Img(\phi_r)=Ae_2r$ which is a nilpotent left ideal in $A$ and since $Ae_2\simeq \nu_A(Ae_1)$ which is indecomposable and injective,  $\phi_r$ is not injective. Otherwise, we would have $Ae_2r\simeq Ae_2$ and $A=Ae_2r\oplus L$ for some left ideal $L$ of $A$ since the module $Ae_2r$ is injective, and consequently $Ae_2r$ would contain an nonzero idempotent element and therefore not be nilpotent, a contradiction. Thus $(\soc(Ae_2))\phi_r=0$, and therefore $ur=0$ and $u\in \soc(e_1A)$. Now, for each $0\neq x\in Ae_2$,  $Ax$ is a nonzero submodule of $Ae_2$.  Hence $\soc(Ae_2)\subseteq Ax$ and $u\in Ax$. Thus there is some $a\in A$ such that $u=ax$. Let $y=e_1a$. Then $u=e_1u=e_1ax=yx$. That is, $u$ satisfies the condition (i). Similarly, one can prove that $u$ satisfies the condition (ii) by the fact that $u\in\soc(e_1A)$ and $e_1A\simeq D(e_2A)$ which is an indecomposable, injective right $A$-module.

$(2)\Rightarrow (3)$ This is trivial.

$(3)\Rightarrow (1)$ Suppose that $u\in e_1\soc (Ae_2)$ is a nonzero element satisfying the conditions (i) and (ii). Let $\alpha$ be a linear map in $D(e_1A)$ such that $(u)\alpha=1$. Define
 $\phi: Ae_2\lra D(e_1A), \quad z\mapsto (z\,\cdot) \alpha.$
 Then $\phi$ is a homomorphism of $A$-modules.
For $0\ne x\in Ae_2$,  there is an element $y\in e_1A$ such that $yx=u$ by the condition (i). It follows that $(x)\phi$ sends $y$ to $(yx)\alpha=(u)\alpha=1$. Thus $(x)\phi\neq 0$. This implies that $\phi$ is injective. Similarly, let $\beta$ be a linear map in $D(Ae_2)$ such that $(u)\beta=1$. Using the condition (ii), one can prove that the map
$$e_1A\lra D(Ae_2),\quad  y\mapsto \beta(\cdot\, y)$$
is an injective  homomorphism of right $A$-modules. Then $\dim_kD(e_1A)=\dim_k e_1A\leq \dim_kD(Ae_2)=\dim_k Ae_2\leq\dim_kD(e_1A).$
It follows that these dimensions are equal, and consequently $\phi$ is an isomorphism.
$\square$

To describe $\nu$-stably projective $\Gamma$-modules, we also need the following lemma.
\begin{Lem}
The assignment $(e)\theta=\overline{e+f}$, $(e_i)\theta=\bar{e}_i$ for all $i\geq 2$, and $(r)\theta=\bar{r}$ for all $r\in r_A$, defines an injective $k$-linear map $\theta: A\ra \Gamma$
such that   $(xy)\theta=(x)\theta (y)\theta$ for all $x, y\in A$.
\label{lemma-map-theta}
\end{Lem}

{\it Proof.}
Since $A=ke\oplus ke_2\oplus\cdots\oplus ke_n\oplus r_A$, the assignment $(e)\phi:=e+f$, $(e_i)\phi:=e_i$ for all $i\geq 2$ and $(r)\phi:=r$ for all $r\in r_A$ defines a $k$-linear map $\phi: A\ra \Lambda$. Now,  it is rather straightforward to check from definition that $\phi$ is injective and satisfies $(xy)\phi=(x)\phi(y)\phi$ for all $x, y\in A$.  Clearly, $\theta$ is the composite $\phi\pi$, where $\pi$ is the canonical surjective algebra homomorphism from $\Lambda$ to $\Gamma$. Note that $\theta$ is  injective since $A\cap I=\{0\}$, and satisfies the other conditions of the lemma.
$\square$

\begin{Prop}
Keep the above notation. We have the following statements:

$(1)$ For each $i\geq 2$, $\Gamma \bar{e}_i$ is $\nu$-stably projective if and only if so is $Ae_i$.

$(2)$ For each $i\geq 2$, $\Gamma \bar{f}_i$ is $\nu$-stably projective if and only if so is $Bf_i$.

$(3)$ $\Gamma(\overline{e+f})$ is $\nu$-stably projective if and only if $\nu_A(Ae)\simeq Ae$ and $\nu_B(Bf)\simeq Bf$.
\label{proposition-gluesoc-nustp}
\end{Prop}

{\it Proof.} We shall frequently use the injective map $\theta: A\ra \Gamma$ in Lemma \ref{lemma-map-theta}.

(1) Since $a\in A$ is a complete $e$-cycle and contained in $e\cdot\soc(Ae)$ by definition, the socle of $Ae$ is isomorphic to the top of $Ae$.  So, it cannot happen that $\nu_A(Ae_i)\simeq Ae$ for any $i\geq 2$.  Thus, if $Ae_i$, with $i\geq 2,$ is $\nu$-stably projective, then $\nu_A^t(Ae_i)$ is isomorphic to some module in $\{Ae_2, \cdots, Ae_n\}$ for all $t\geq 1$.  Similarly, $\bar{a}$ is a complete $\overline{e+f}$-cycle in $\Gamma$, and it cannot happen that $\nu_{\Gamma}(\Gamma\bar{e}_i)\simeq \Gamma(\overline{e+f})$.  It is also impossible that $\nu_{\Gamma}(\Gamma\bar{e}_i)\simeq \Gamma\bar{f}_l$ for any $l\geq 2$, since $\bar{e}_i\Gamma\bar{f}_l=0$.
We shall show, for $i, j\geq 2$, that   $\nu_A(Ae_i)\simeq Ae_j$  if and only if $\nu_{\Gamma}(\Gamma\bar{e}_i)\simeq\Gamma\bar{e}_j$.  Note that $\theta$ induces isomorphisms of vector spaces $e_iA\ra \bar{e}_i\Gamma$ and $Ae_j\ra \Gamma\bar{e}_j$. Then it is easy to check that a nonzero element $u\in e_iAe_j$ satisfies both (i) and (ii) in Lemma \ref{lemma-Ae-Nakayama} if and only if $(u)\theta$, which is a nonzero element in $\bar{e}_i\Gamma\bar{e}_j$, satisfies the same conditions. By Lemma \ref{lemma-Ae-Nakayama}, $\nu_A(Ae_i)\simeq Ae_j$ if and only if $\nu_{\Gamma}(\Gamma\bar{e}_i)\simeq \Gamma\bar{e}_j$.  Repeating this process,  we see that $\Gamma\bar{e}_i$ is $\nu$-stably projective if and only if
$Ae_i$ is $\nu$-stably projective. This proves (1).

(2) This can be shown similarly.

(3) We assume that $\Gamma(\overline{e+f})$ is $\nu$-stably projective. Then $\Gamma(\overline{e+f})$ is indecomposable and projective-injective, and has a $1$-dimensional simple socle. The element $\bar{a}=\bar{b}$ is in the socle of $\Gamma(\overline{e+f})$ and $(\overline{e+f})\bar{a}=\bar{a}$. It follows that the socle of $\Gamma(\overline{e+f})$ is isomorphic to the simple $\Gamma$-module corresponding to the primitive idempotent $\overline{e+f}$. Hence $\nu_{\Gamma}\Gamma(\overline{e+f})\simeq \Gamma(\overline{e+f})$. By Lemma \ref{lemma-Ae-Nakayama}, the element $\bar{a}\in (\overline{e+f})\Gamma(\overline{e+f})$ satisfies the condition (i) and (ii) in Lemma \ref{lemma-Ae-Nakayama}(3). We shall prove that $a\in eAe$ satisfies the conditions (i) and (ii) in Lemma \ref{lemma-Ae-Nakayama}(3).
Let $x\in Ae$ be a  nonzero element.  Suppose $x=\lambda e+r$ for some $\lambda\in k$ and $r\in r_Ae$.  If $\lambda\neq 0$, then $\frac{1}{\lambda}ax=a+\frac{1}{\lambda}ar=a$.
Now, we assume that $\lambda=0$ and $x=r\in r_Ae$. In this case,  $(x)\theta=\bar{r}$ is a non-zero element in $\Gamma(\overline{e+f})$. By Lemma \ref{lemma-Ae-Nakayama}, there is an element $w=\mu(\overline{e+f})+\bar{r}_1+\bar{r}_2$ in $(\overline{e+f})\Gamma$, where $\mu\in k$, $r_1\in er_A$ and $r_2\in fr_B$, such that $w\cdot (x)\theta=\bar{a}$. Note that $\bar{r}_2\cdot (x)\theta=\bar{r}_2\bar{r}=0$. Hence $\bar{a}=w\cdot (x)\theta=(\mu(\overline{e+f})+\bar{r}_1)\cdot (x)\theta=(\mu e+r_1)\theta\cdot (x)\theta$. Let $y=\mu e+r_1$. Then $y\in eA$ and $(yx)\theta=\bar{a}$. Since $a\in r_A$, we have $\bar{a}=(a)\theta$. It follows that $(yx)\theta=(a)\theta$, and consequently $yx=a$. This shows that $a$ satisfies the condition (i) in Lemma \ref{lemma-Ae-Nakayama}. Similarly, $a$ also satisfies the condition (ii) in Lemma \ref{lemma-Ae-Nakayama}. Consequently, $\nu_A(Ae)\simeq Ae$.  Similarly, $\nu_B(Bf)\simeq Bf$.

Conversely, we assume that $\nu_A(Ae)\simeq Ae$ and $\nu_B(Bf)\simeq Bf$. By Lemma \ref{lemma-Ae-Nakayama}, $a\in eAe$ (respectively, $b\in fBf$) satisfies the conditions (i) and (ii) in Lemma \ref{lemma-Ae-Nakayama}.  We claim that $\bar{a}\in (\overline{e+f})\Gamma(\overline{e+f})$ also satisfies the conditions (i) and (ii) in Lemma \ref{lemma-Ae-Nakayama}. Let $w=\lambda(\overline{e+f})+\bar{r}_1+\bar{r}_2\in \Gamma(\overline{e+f})$ be a nonzero element with $\lambda\in k$, $r_1\in r_Ae$ and $r_2\in r_Bf$. If  $\lambda\neq 0$, $\frac{1}{\lambda}\bar{a}\in (\overline{e+f})\Gamma$, and $\frac{1}{\lambda}\bar{a}w=\bar{a}$. Next we assume $\lambda=0$.  Then $w=\bar{r}_1+\bar{r}_2$.  We can assume that $w$ is not a multiple of $\bar{a}=\bar{b}$ (Otherwise, $w$ clearly satisfies the condition (i) in Lemma \ref{lemma-Ae-Nakayama}). Then either $\bar{r}_1$ or $\bar{r}_2$ is not a multiple of $\bar{a}$. Without loss of generality, we assume that $\bar{r}_1\neq 0$ is not a multiple of $\bar{a}$.  This is equivalent to $Ar_1\neq ka$. Since $ka$ is the simple socle of $Ae$ and $Ar_1$ is a nonzero submodule of $Ae$, $ka$ is a proper submodule of $Ar_1$, and there is some $z\in r_A$ such that $zr_1=a$. Let $y=ez$. Then $yr_1=ezr_1=ea=a$. Hence $\bar{y}\cdot w=\bar{y}\bar{r}_1+\bar{y}\bar{r}_2=\bar{a}$.  Altogether, $\bar{a}$ satisfies the condition (i) in Lemma \ref{lemma-Ae-Nakayama}. Similarly, $\bar{a}$ also satisfies the condition (ii) in Lemma \ref{lemma-Ae-Nakayama}. As a result, $\nu_{\Gamma}\big(\Gamma(\overline{e+f})\big)\simeq \Gamma(\overline{e+f})$.
$\square$

Note that if both $A$ and $B$ are symmetric algebras, then Proposition \ref{proposition-gluesoc-nustp} tells that $A\gluesocab B$ is a self-injective algebra. Further,
we have the following corollary of Proposition \ref{proposition-gluesoc-nustp} and Corollary \ref{cor-identifysoc}. For notation, see Section \ref{subsection-glue-e-cycles}.

\begin{Koro}
Let $A$ be a Frobenius-free $k$-algebra given by a quiver with relations. Suppose that $e\in A$ is a primitive idempotent and that $a\in A$ is a complete $e$-cycle element. Let $\cpx{T}$ be a basic, radical tilting complex over $A$ such that $[\cpx{T}:Ae]=1$, that $A':=\End_{\Kb{A}}(\cpx{T})$ is Frobenius-finite, and that the only $\nu$-stably projective, indecomposable $A'$-module is
$A'\tilde{e}$ with $ \tilde{e}^2 = \tilde{e}\in A'$, where $\tilde{e}$ is the idempotent element in $A'$ corresponding to the indecomposable direct summand of $\cpx{T}$ in which $Ae$ appears. Let $B$ be a basic, symmetric $k$-algebra which is splitting over $k$ and has no nonzero semisimple direct summands, and let $f$ be a primitive idempotent in $B$ and $0\ne b\in \soc(_BBf)$.  Then $A\gluesocab B$  and $A'\gluesoc{a'}{b}B$ are derived equivalent, and their Frobenius parts are $(1_B-f)B(1_B-f)$ and $B$, respectively. \label{cor3}
\end{Koro}

Finally, we employ Corollary \ref{cor3} to construct a series of examples, showing that Frobenius type of algebras may change under derived equivalences.
\begin{Bsp}{\rm Let $A$ and $A'$ be $k$-algebras given by the following quivers with relations, respectively:}\label{Bsp5.6}
\end{Bsp}
\vspace{-0.3cm}
$$\begin{array}{ccc}
\xymatrix{
 \bullet\ar@<2.5pt>[r]^{\alpha}^(0.1){1}^(1){2}& \bullet\ar@<2.5pt>[l]^{\delta}\ar@<2.5pt>[r]^{\beta}^(1){3} &\bullet\ar@<2.5pt>[l]^{\gamma}}
 &\hspace{2cm}&
\xymatrix{
\bullet\ar[r]^{\alpha'}^(0){1}^(1){2}& \bullet\ar[dl]^{\beta'}\\
\bullet\ar[u]^{\gamma'}^(0){3} }\\
\alpha\delta\alpha, \gamma\delta, \delta\alpha-\beta\gamma\; ;
&&
 \alpha'\beta'\gamma'\alpha',\, \gamma'\alpha'\beta'\gamma'. \\ \end{array}$$
Then $\dim_k(A)=12$ and $\dim_k(A')=13.$ We denote by $e_i$ the primitive idempotent element of $A$ corresponding to the vertex $i$.  Let $e=e_1$.  Then there is a tilting complex $\cpx{T}=\cpx{T_1}\oplus Ae_2[1]\oplus Ae_3[1]$ over $A$, where $\cpx{T_1}$ is the complex
 $$0\lra Ae_2\lraf{\cdot \delta} Ae\lra 0$$
 with $Ae$ in degree zero. Then the assignment
 $$ \begin{array}{ccc}
 \alpha'\quad\mapsto\quad \xymatrix@M=0.6mm@C=5mm{
 0\ar[r] &Ae_2\ar@{=}[d]\ar[r]^{\cdot \delta} & Ae\ar[d]\ar[r]& 0\\
0\ar[r] &Ae_2\ar[r] &0
 } & \quad &  \beta'\quad\mapsto\quad \xymatrix@M=0.6mm@C=5mm{
 0\ar[r] &Ae_2\ar[d]^{\cdot \beta}\ar[r]& 0\\
0\ar[r] &Ae_3\ar[r] &0
 }  \\
 &&\\
 \gamma'\quad\mapsto\quad \xymatrix@M=0.6mm@C=5mm{
 0\ar[r] &Ae_3\ar[d]^{\cdot (-\gamma)}\ar[r] &0\ar[d]\\
 0\ar[r] &Ae_2\ar[r]^{\cdot \delta} & Ae\ar[r]& 0
 } & \quad & \\
 \end{array}$$
induces an isomorphism between $A'$ and $\widetilde{A}:=\End_{\Kb{A}}(\cpx{T})$. In the following, we identify $A'$ with $\End_{\Kb{A}}(\cpx{T})$, that is, $A' = \widetilde{A}$. The element $a:=\alpha\delta$ is a complete $e$-cycle in $A$.  Let $a'$ be the particular element in $\widetilde{A}$:
$$\xymatrix@M=0.6mm{
0\ar[r]  & Ae_2\ar[d]^{0}\ar[r]^{\cdot \delta} &Ae\ar[d]^{\cdot \alpha\delta}\ar[r] &0\\
0\ar[r] & Ae_2\ar[r]^{\cdot \delta} &Ae\ar[r] &0.\\
}$$
Due to the relation $\delta\alpha-\beta\gamma$ in $A$, we get a commutative diagram
$$\xymatrix@M=0.5mm{
0\ar[r]  & Ae_2\ar[d]_{\cdot \beta\gamma}\ar[r]^{\cdot \delta} &Ae\ar@{-->}[ld]^{\cdot \alpha}\ar[d]^{\cdot \alpha\delta}\ar[r] &0\\
0\ar[r] & Ae_2\ar[r]_{\cdot \delta} &Ae\ar[r] &0.\\
}$$
This implies $a'=\alpha'\beta'\gamma'$.

Let $\tilde{e}$ be the primitive idempotent in $\widetilde{A}$ corresponding to the direct summand $\cpx{T_1}$. Under the identification of $\tilde{A}$ with $A'$, $\tilde{e}$ corresponds to the vertex $1$ in $A'$. Note that $A$ is Frobenius-free, the Frobenius part of $A'$ is isomorphic to $k[x]/(x^2)$, and the only $\nu$-stably projective, indecomposable $A'$-module is $A'\tilde{e}$, that is $A'$-stp $=\add(A'\tilde{e})$.

Let $B$ be a basic, symmetric algebra without semisimple direct summands, and let $f$ be a primitive idempotent in $B$.  Then any nonzero element in the socle of $Bf$ is a complete $f$-cycle.  Let $b$ be such an element. Then, by Corollary \ref{cor3}, the algebras $A\gluesocab B$ and $A'\gluesoc{a'}{b}B$ are derived equivalent, while the Frobenius part of $A\gluesocab B$ is $(1_B-f)B(1_B-f)$ and the Frobenius part of $A'\gluesoc{a'}{b}B$ is $B$.

Thus, if we choose a basic, symmetric algebra $B$ and a primitive idempotent $f\in B$ such that $B$ is wild (or tame) and $(1_B-f)B(1_B-f)$ is tame or representation-finite, then $A'\gluesoc{a'}{b}B$ is Frobenius-wild (or tame), while $A\gluesocab B$ is Frobenius-tame or Frobenius-finite.
This means that Frobenius type may change under derived equivalences in general.

Note that under the derived equivalence defined by $\cpx{T}$ in Example \ref{Bsp5.6}, the simple $A$-module $S_1$ corresponding to the vertex $1$ is sent to the simple $A'$-module $S'_1$ corresponding to the vertex $1$ by Lemma \ref{LemsimpleTosimple}. Thus, by Theorem \ref{TheoremGluing}, if we glue an arbitrary algebra $B$ at the vertex $1$ in $A$ and $A'$, respectively, then the resulting extension algebras of $A$ and $A'$ are also derived equivalent.

\smallskip
{\footnotesize
\bigskip Wei Hu, School of Mathematical Sciences, Beijing Normal
University, 100875 Beijing, People's Republic of  China

{\tt Email: huwei@bnu.edu.cn}

\bigskip
Changchang Xi, School of Mathematical Sciences, Capital Normal
University, 100048 Beijing, People's Republic of  China

{\tt Email: xicc@cnu.edu.cn}}
\end{document}